\documentclass[english,11pt]{article}

\usepackage[german,french,english]{babel}
\usepackage[cp850]{inputenc}
\usepackage{latexsym,graphicx, fancybox}
\usepackage{amssymb, amsmath, amsfonts}
\usepackage{color}
\usepackage{latexsym}

\font\tenmath=msbm10 scaled 1200
\font\sevenmath=msbm7 scaled 1200
\font\fivemath=msbm5 scaled 1200 
\newfam\mathfam \textfont\mathfam=\tenmath
\scriptfont\mathfam=\sevenmath \scriptscriptfont\mathfam=\fivemath

\def\R{{\mathbb R}}
\def\N{{\mathbb N}}
\def\E{{\mathbb E}}
\def\L{{\cal L}}
\def\P{{\mathbb P}}

\newtheorem{Theorem}{Theorem}[section]

\newtheorem{Proposition}[Theorem]{Proposition}

\newtheorem{Lemma}[Theorem]{Lemma}
\newtheorem{Corollary}[Theorem]{Corollary}
\newtheorem{Remark}[Theorem]{Remark}
\newtheorem{Example}[Theorem]{Example}

\newfam\mathfam \textfont\mathfam=\tenmath
\scriptfont\mathfam=\sevenmath \scriptscriptfont\mathfam=\fivemath

\def \^#1{\if#1i{\accent"5E\i}\else{\accent"5E#1}\fi}

\def \ind {1 \mkern -5mu \hbox{I}}

\def \b{\beta}

\begin{document}
\selectlanguage{english}
\title{\bf Convex comparison of Gaussian mixtures}
 
\author{ 
{\sc Benjamin Jourdain} \thanks{CERMICS, Ecole des Ponts, IP Paris, INRIA, Marne-la-Vall\'ee, France. E-mail: {\tt   benjamin.jourdain@enpc.fr}}~\footnotemark[3]
\and   
{\sc  Gilles Pag\`es} \thanks{Laboratoire de Probabilit\'es, Statistique et Mod\'elisation, UMR~8001, Campus Pierre et Marie Curie, Sorbonne Universit\'e case 158, 4, pl. Jussieu, F-75252 Paris Cedex 5, France. E-mail: {\tt  gilles.pages@upmc.fr}}~\thanks{This research
benefited from the support of the ``Chaire Risques Financiers'', Fondation du Risque}  }
\maketitle
\begin{abstract}
   Motivated by the study of the propagation of convexity by semi-groups of stochastic differential equations and convex comparison between the distributions of solutions of two such equations, we study the comparison for the convex order between a Gaussian distribution and a Gaussian mixture. We give and discuss intrinsic necessary and sufficient conditions for convex ordering. On the examples that we have worked out, the two conditions appear to be closely related.
\end{abstract}\section*{Introduction}
A successfull strategy implemented in \cite{BRFS,BRAP,PSP,LiuPag,LiuPagAAP,JourPagVolt22,JouPaSDE} in order to establish propagation of convexity by solutions of stochastic differential equations or stochastic volterra equations and convex comparison between the solutions of two such equations is to first deal with Euler discretizations of these equations and then take the limit as the time-step goes to $0$. To deal with the Euler schemes, one needs conditions ensuring that for some diffusion coefficient $\sigma:\R^d\to\R^{d\times q}$ and a Gaussian random vector $G\sim{\cal N}_q(0,I_q)$, the mapping which associates to $x\in\R^d$ the law ${\cal L}(\sigma(x)G)$ of the random vector $\sigma(x)G$ is convex when the codomain is endowed with the convex order. For $\xi\in\R^d$ and $\Sigma\in\R^{d\times d}$ symmetric positive semi-definite, ${\cal N}_d(\xi,\Sigma)$ denotes the $d$-dimensional Gaussian law with expectation $\xi$ and covariance matrix $\Sigma$ and the convex order between two probability measures $\mu$ and $\nu$ on $\R^d$ with finite first-order moment is denoted by $\mu\le_{cx}\nu$ and defined by \begin{equation}
   \forall \varphi:\R^d\to\R\mbox{ convex },\;\int_{\R^d}\varphi(x)\mu(dx)\le\int_{\R^d}\varphi(x)\nu(dx).\label{defocv}
\end{equation}
In \cite{JP}, we show that a sufficient condition for the
inequality $${\cal L}(\sigma(p_1x+(1-p_1)y)Z)\le_{cx}p_1{\cal L}(\sigma(x)Z)+(1-p_1){\cal L}(\sigma(y)Z)$$
where $p_1\in (0,1)$, $x,y\in\R^d$ and $Z$ is any integrable $\R^q$-valued random vector with radial distribution (i.e. such that ${\cal L}(OZ)={\cal L}(Z)$ for each orthogonal matrix $O\in\R^{q\times q}$, e.g. $Z=G$) is 
\begin{align}
   &\exists O\in\R^{q\times q}\mbox{ orthogonal such that }\notag\\&\sigma\sigma^*(p_1x+(1-p_1)y)\le (p_1\sigma(x)+(1-p_1)\sigma(y)O)(p_1\sigma(x)+(1-p_1)\sigma(y)O)^*\label{condconvjp}
\end{align}
where the inequality means that the difference between the right-hand and the left-hand sides is a positive semi-definite matrix.
This generalized the sufficient condition in \cite{LiuPagAAP} where $O$ was chosen equal to the identity matrix $I_q\in\R^{q\times q}$.

In the Gaussian noise case $Z=G$, denoting $\Sigma=\sigma\sigma^*(p_1x+(1-p_1)y)$,  $\Sigma_1=\sigma\sigma^*(x)$ and $\Sigma_2=\sigma\sigma^*(y)$, one immediately sees that intrisic conditions for the inequality $${\cal N}_d(0,\Sigma)\le_{cx}p_1{\cal N}_d(0,\Sigma_1)+(1-p_1){\cal N}_d(0,\Sigma_2)$$ should only involve $p_1,\Sigma,\Sigma_1$ and $\Sigma_2$ and not the matrices $\sigma(p_1x+(1-p_1)y),\sigma(x)$ and $\sigma(y)$.
In the present paper, we investigate and propose such an intrinsic sufficient condition and an intrinsic necessary condition for the more general inequality
$${\cal N}_d(0,\Sigma)\le_{cx}\sum_{i=1}^n p_i{\cal N}_d(0,\Sigma_i),$$
where $n\ge 2$, $(p_1,\cdots,p_n)\in(0,1)^n$ satisfy $\sum_{i=1}^n p_i=1$ and $\Sigma,\Sigma_1,\cdots,\Sigma_n\in\R^{d\times d}$ are covariance matrices i.e. positive semi-definite matrices. 
Our sufficient condition is the existence of a $nd$-dimensional Gaussian coupling between random vectors following the distributions in the mixture such that the law of the convex combination with weights $(p_1,\cdots,p_n)$ of these random vectors dominates ${\cal N}_d(0,\Sigma)$ for the convex order. The necessary condition relies on \eqref{defocv} written for the convex functions given by the absolute value of linear functions $\R^d\ni x\mapsto\varphi(x)=|\xi^*x|$ with $\xi\in\R^d$. It turns out that these conditions keep sufficient/necessary for
\begin{equation}
   \forall x_1,\cdots,x_n\in\R^d\mbox{ s.t. }\sum_{i=1}^n p_ix_i=0,\;  {\cal N}_d(0,\Sigma)\le_{cx}\sum_{i=1}^np_i{\cal N}_d(x_i,\Sigma_i).\label{cc}
\end{equation}
Even if we have not been able to exhibit an example where the necessary condition is satisfied but the sufficient condition is not, we believe that the latter is stronger than the former. On the other hand, all the situations where we could check equivalence between the two conditions and therefore with \eqref{cc} suppose that the same correlation matrix is associated with $M\Sigma_1M^*,\cdots,M\Sigma_nM^*$ for some non-singular matrix $M\in\R^{d\times d}$.

In the first section, we prove and discuss these conditions. The second section is dedicated to the case where there are only $n=2$ Gaussian distributions in the mixture which was our initial motivation. We can make the conditions more explicit in this particular setting.

\noindent {\sc Notation}

\smallskip
\noindent 
$\bullet$ We denote $\R^{r\times c}$ the set of real matrices with $r$ rows and $c$ columns. 

\smallskip
\noindent $\bullet$ For $A\in \R^{r\times c}$, we denote  $A^*\in \R^{c\times r}$ the transpose matrix of $A$.

\smallskip
\noindent $\bullet$ For $A\in \R^{r\times c}$, $1\le r_1\le r_2\le r$ and $1\le c_1\le c_2\le c$ we denote by $A_{r_1:r_2,c_1:c_2}\in\R^{(r_2-r_1+1)\times (c_2-c_1+1)}$ the submatrix whose entry with index $ij$ is equal to $A_{(r_1+i-1)(c_1+j-1)}$.

\smallskip
\noindent $\bullet$ We denote by ${\cal I}(d)=\{A\in\R^{d\times d}\mbox{ non singular}\}$ the set of matrices $A$ which admit an inverse $A^{-1}$.

\smallskip
\noindent $\bullet$ We denote by ${\cal S}_+(d)= \big\{S\in \R^{d\times d}, \mbox{ symmetric, positive semi-definite}\big\}$ the set of covariance matrices, by ${\cal C}(d)= \big\{C\in {\cal S}_+(d):C_{ii}=1\mbox{ for }\in\{1,\cdots,d\}\big\}$ the subset of correlation matrices and by ${\cal O}(d)= \big\{ O\in \R^{d\times d}:  OO^* =I_d\big\}$ where $I_d\in{\cal S}_+(d)$ is the identity matrix, the set of orthogonal matrices.

\smallskip
\noindent $\bullet$ For $\lambda_1,\cdots,\lambda_d\in\R^d$, we denote by ${\rm diag}(\lambda_1\cdots,\lambda_d)$ the diagonal matrix in $\R^{d\times d}$ with successive diagonal coefficients $\lambda_1,\cdots,\lambda_d$.

\smallskip
\noindent $\bullet$  For $S,T\in \R^{d\times d}$ symmetric, $S\le T$ means that $T-S\in {\cal S}_+(d)$.

\smallskip
\noindent $\bullet$  For $S\in {\cal S}_+(d)$ (resp. $S\in {\cal S}_+(d)$  non sigular), we denote by ${S}^{1/2}$ (resp. $S^{-1/2}$) the unique element of ${\cal S}_+(d)$ such that ${S}^{1/2}{S}^{1/2}=S$ (resp. $S^{-1/2}S^{1/2}=I_d$).

\section{Convex ordering between a Gaussian distribution and a Gaussian mixture}\label{sec:discutconv}
\subsection{Main results and comments}
For $\Sigma,\Sigma_1,\cdots,\Sigma_n\in{\cal S}_+(d)$, we want a condition ensuring that \begin{equation*}
   {\cal N}_d(0,\Sigma)\le_{cx}\sum_{i=1}^np_i{\cal N}_d(0,\Sigma_i).
\end{equation*}
We recall that a correlation matrix associated with a matrix $\Theta\in{\cal S}_+(d)$ is a matrix $C\in{\cal C}(d)$ such that
$$\Theta={\rm diag}\left(\sqrt{\Theta_{11}},\cdots,\sqrt{\Theta_{dd}}\right)C{\rm diag}\left(\sqrt{\Theta_{11}},\cdots,\sqrt{\Theta_{dd}}\right).$$
The matrix $C$ with all diagonal entries equal to $1$ and all non-diagonal entries equal to those of
${\rm diag}\left(\frac{1_{\{\Theta_{11}>0\}}}{\sqrt{\Theta_{11}}},\cdots,\frac {1_{\{\Theta_{dd}>0\}}}{\sqrt{\Theta_{dd}}}\right)\Theta{\rm diag}\left(\frac{1_{\{\Theta_{11}>0\}}}{\sqrt{\Theta_{11}}},\cdots,\frac {1_{\{\Theta_{dd}>0\}}}{\sqrt{\Theta_{dd}}}\right)$ is a correlation matrix associated with $\Theta$. When $\Theta$ is non singular, it is equal to ${\rm diag}\left(\frac 1{\sqrt{\Theta_{11}}},\cdots,\frac 1{\sqrt{\Theta_{dd}}}\right)\Theta{\rm diag}\left(\frac 1{\sqrt{\Theta_{11}}},\cdots,\frac 1{\sqrt{\Theta_{dd}}}\right)$ and this is the only correlation matrix associated with $\Theta$. 
\begin{Theorem}\label{proprinc}
  Let $n\ge 2$, $\Sigma,\Sigma_1,\cdots,\Sigma_n\in{\cal S}_+(d)$ and $(p_1,\cdots,p_n)\in(0,1)^n$ such that $\sum_{i=1}^n p_i=1$. Then \begin{align}
  &\exists (M,C)\in{\cal I}(d)\times{\cal C}(d)\mbox{ such that }C\mbox{ is associated with }M\Sigma_iM^*\mbox{ for }i\in\{1,\cdots,n\}\mbox{ and }\notag\\& M\Sigma M^*\le DCD\mbox{ where }D={\rm diag}\left(\sum_{i=1}^np_i\sqrt{(M\Sigma_i M^*)_{11}},\cdots,\sum_{i=1}^np_i\sqrt{(M\Sigma_i M^*)_{dd}}\right)\label{correl}
\end{align}implies 
\begin{align}
 \exists\Gamma\in{\cal S}_+(nd)&\mbox{ such that }\Sigma\le (p_1I_d,\cdots,p_n I_d)\Gamma (p_1I_d,\cdots,p_n I_d)^*\notag\\&
  \mbox{ and }\forall i\in\{1,\cdots,n\},\;\Gamma_{(i-1)d+1:id,(i-1)d+1:id}=\Sigma_i\label{inecov}
\end{align}
which implies \begin{equation}
 \forall x_1,\cdots,x_n\in\R^d\mbox{ s.t. }\sum_{i=1}^n p_ix_i=0,\;  {\cal N}_d(0,\Sigma)\le_{cx}\sum_{i=1}^np_i{\cal N}_d(x_i,\Sigma_i)\label{compccgg}
 \end{equation} which in turn implies \begin{equation}
   \forall \xi\in\R^d,\;\sqrt{\xi^*\Sigma\xi}\le \sum_{i=1}^n p_i\sqrt{\xi^*\Sigma_i\xi}.\label{inegsqrt}
 \end{equation}
Moreover, the four conditions \eqref{correl}, \eqref{inecov}, \eqref{compccgg} and \eqref{inegsqrt} are equivalent \begin{itemize}
\item when the matrices $\Sigma_i,\;i\in\{1,\cdots,n\}$ are all colinear and in particular in dimension $d=1$,
\item  or when there exists $N\in{\cal I}(d)$ such that $N\Sigma_iN^*N\Sigma_jN^*=0$ ($\Leftrightarrow \Sigma_iN^*N\Sigma_j=0$) for $i,j\in\{1,\cdots,n\}$ such that $i\neq j$ and in particular when $\Sigma_i\Sigma_j=0$ for $i,j\in\{1,\cdots,n\}$ such that $i\neq j$,
  \item or % \textcolor{red}{when there exists $M\in\R^{d\times d}$ non singular such that $M\Sigma M^*,M\Sigma_1M^*,\cdots,M\Sigma_nM^*$ all are diagonal} and in particular
    when there exists $N\in{\cal I}(d)$ such that the matrices $N\Sigma N^*,N\Sigma_1N^*,\cdots,N\Sigma_nN^*$ commute and in particular when $\Sigma,\Sigma_1,\cdots,\Sigma_n$ commute or when there exists $M\in{\cal I}(d)$ such that the matrices $M\Sigma M^*,M\Sigma_1M^*,\cdots,M\Sigma_nM^*$ are diagonal,
  \item or when there exists $(M,C)\in{\cal I}(d)\times{\cal C}(d)$ such that $C$ is associated with all the matrices $M\Sigma_1 M^*,\cdots, M\Sigma_nM^*$ and with the matrix $\widehat \Sigma\in\R^{d\times d}$ with diagonal entries $\widehat \Sigma_{kk}=\left(\sum_{i=1}^np_i\sqrt{(M\Sigma_i M^*)_{kk}}\right)^2$ for $k\in\{1,\cdots,d\}$ and all other entries equal to those of $M\Sigma M^*$.
\end{itemize} 
\end{Theorem}
As an obvious consequence of the equivalence in dimension $d=1$, we deduce that
\begin{Corollary}\label{cor1d}Let $\sigma,\sigma_1,\cdots,\sigma_n\in\R_+$ and $(p_1,\cdots,p_n)\in(0,1)^n$ be such that $\sum_{i=1}^n p_i=1$. Then
  $$\forall x_1,\cdots,x_n\in\R\mbox{ s.t. }\sum_{i=1}^n p_ix_i=0,\;{\cal N}_1(0,\sigma^2)\le_{cx} \sum_{i=1}^n p_i{\cal N}_1(x_i,\sigma^2_i)\Leftrightarrow \sigma\le \sum_{i=1}^n p_i\sigma_i.$$
\end{Corollary}
The proof of Theorem \ref{proprinc} is given in Section \ref{sec:proof}.
When considering the converse inequality where a Gaussian mixture is dominated for the convex order by a single Gaussian distribution, we now derive a simple necessary and sufficient condition, also proved in Section \ref{sec:proof}.\begin{Proposition}\label{propcnsconvcomble}
 Let $n\ge 2$, $\Sigma,\Sigma_1,\cdots,\Sigma_n\in{\cal S}_+(d)$ and $(p_1,\cdots,p_n)\in(0,1)^n$ be such that $\sum_{i=1}^n p_i=1$. Then
  \begin{align*}
    \sum_{i=1}^np_i{\cal N}_d(0,\Sigma_i)\le_{cx}{\cal N}_d(0,\Sigma)\Leftrightarrow \forall i\in\{1,\cdots,n\},\;\Sigma_i\le \Sigma.\end{align*}
\end{Proposition}
\begin{Remark}
   When $\Sigma\le \Sigma_i$ for all $i\in\{1,\cdots,n\}$, then \eqref{inecov} holds with $\Gamma$ such that $$\Gamma_{(i-1)d+1:id,(j-1)d+1:jd}=1_{\{i=j\}}\Sigma_i+1_{\{i\ne j\}}\Sigma\mbox{ for }i,j\in\{1,\cdots,n\}.$$ This matrix belongs to ${\cal S}_+(nd)$ since it is larger than $(I_d,\cdots,I_d)^*\Sigma(I_d,\cdots,I_d)$.
\end{Remark}
  
\begin{Remark}
  For $X_i\sim{\cal N}_d(0,\Sigma_i)$ and $\varphi:\R^d\to\R$ convex and even, $\inf_{x_i\in\R^d}\E[\varphi(X_i+x_i)]=\E[\varphi(X_i)]$ since $$\varphi(X_i)\le\frac 12\left(\varphi(X_i+x_i)+\varphi(X_i-x_i)\right)=\frac 12\left(\varphi(X_i+x_i)+\varphi(-X_i+x_i)\right)$$ and  $-X_i\sim{\cal N}_d(0,\Sigma_i)$. This explains why, in the proof of \eqref{compccgg}$\Rightarrow$\eqref{inegsqrt} which relies on the choice $\varphi(x)=|\xi^*x|$, we only use ${\cal N}_d(0,\Sigma)\le_{cx}\sum_{i=1}^np_i{\cal N}_d(0,\Sigma_i)$.

  However, the inequality $\inf_{x_i\in\R^d}\E[\varphi(X_i+x_i)]=\E[\varphi(X_i)]$ does not generalize without evenness of the convex function $\varphi$. Indeed such a generalization would imply that for $x_1,\cdots,x_n\in\R^d\mbox{ s.t. }\sum_{i=1}^n p_ix_i=0$, $\sum_{i=1}^np_i{\cal N}_d(0,\Sigma_i)\le_{cx}\sum_{i=1}^np_i{\cal N}_d(x_i,\Sigma_i)$ while for $\sigma_1,\sigma_2\ge 0$, $p_1\in(0,1)$ and $\lambda\in \R$, the function which associates to $x_1\in\R$ the integral $p_1e^{\lambda x_1+\frac{\lambda^2\sigma_1^2}{2}}+(1-p_1)e^{-\frac{\lambda p_1x_1}{1-p_1}+\frac{\lambda^2\sigma_2^2}{2}}$ of the convex function $\R\ni x\mapsto e^{\lambda x}$ with respect to $p_1{\cal N}_1(x_1,\sigma^2_1)+(1-p_1){\cal N}_1(-\frac{p_1x_1}{1-p_1},\sigma^2_2)$ is minimal for $x_1=\frac 1 2\lambda(1-p_1)(\sigma_2^2-\sigma_1^2)$. 
\end{Remark}
\begin{Remark}
   In Example \ref{exdiag2d} below, we will check that for $\Sigma=\left(\begin{array}{cc}a & b
   \\ b & a
\end{array}\right)$ with $|b|\le a< 3-2\sqrt{2}$, $\Sigma_1=\left(\begin{array}{cc}8 & 0
   \\ 0 & 4
\end{array}\right)$, $\Sigma_2=\left(\begin{array}{cc}4 & 0
   \\ 0 & 8
                                     \end{array}\right)$ and $p_1=p_2=\frac 12$, \eqref{inecov} is satisfied. Since no choice of $c\in\R$ ensures that  the matrix $\left(\begin{array}{ccc} 8 & 0 &2(a-3)\\0 & 4 & c\\2(a-3) & c  &4\end{array}\right)$ with determinant $16(3+2\sqrt{2}-a)(a+2\sqrt{2}-3)-8c^2$ is positive semi-definite, there is no $\Gamma\in{\cal S}_+(4)$ such that $\Gamma_{1:2,1:2}=\Sigma_1$, $\Gamma_{3:4,3:4}=\Sigma_2$ and $\Sigma=(p_1I_2,p_2I_2)\Gamma (p_1I_2,p_2I_2)^*$. Therefore replacing the inequality in \eqref{inecov} by an equality leads to a stronger condition.
\end{Remark}
\begin{Remark}\label{remskew}
  Since $\forall \alpha>0,\;\sqrt{\xi^*\Sigma_i\xi\times\xi^*\Sigma_j\xi}\le \frac{\alpha}2  \xi^*\Sigma_i\xi+\frac 1{2\alpha}\xi^*\Sigma_j\xi$ with equality attained for $\alpha=\sqrt{\frac{\xi^*\Sigma_j\xi}{\xi^*\Sigma_i\xi}}$ when $\xi^*\Sigma_i\xi>0$ and in the limit $\alpha\to\infty$ when $\xi^*\Sigma_i\xi=0$,
  \eqref{inegsqrt} also writes
\begin{equation}
   \forall M\in\R^{n\times n}\mbox{ skew-symmetric},\;\Sigma\le \sum_{i=1}^n p_i\sum_{j=1}^n e^{M_{ij}}p_j\Sigma_i.\label{inegsqrtskew}
\end{equation}% \textcolor{red}{Condition \eqref{inegsqrt} also writes
  % $$\forall \alpha\in (0,+\infty)^{n\times n},\;\Sigma\le \sum_{i=1}^n p_i\sum_{j=1}^n \left(\frac{\alpha_{ij}}{2}+\frac 1{2\alpha_{ji}}\right)p_j\Sigma_i,$$ où on peut prendre $\alpha_{ij}=\frac 1{\alpha_{ji}}$ et en particulier $\alpha_{ii}=1$.}
\end{Remark}
\begin{Remark}\label{reminegsqrt}
   In view of Corollary \ref{cor1d}, the inequality \eqref{inegsqrt} is equivalent to
   $\forall \xi\in\R^d,\;\forall y_1,\cdots,y_n\in\R\mbox{ s.t. }\sum_{i=1}^n p_iy_i=0,\;{\cal N}_1(0,\xi^*\Sigma\xi)\le \sum_{i=1}^n p_i{\cal N}_1(y_i,\xi^*\Sigma_i\xi)$. Note that any convex quadratic function on $\R^d$ writes $\varphi(x)=x^*Mx+\xi_0^*x+c$ for some constant $c\in\R$, some vector $\xi_0\in \R^d$, and some matrix $M\in{\cal S}_+(d)$. By diagonalization of this matrix, $\varphi(x)=\sum_{k=1}^d\lambda_k(\xi_k^*x)^2+\xi_0^*x+c$ for some orthonormal basis $(\xi_k)_{1\le k\le d}$ of $\R^d$ and $\lambda_1,\cdots,\lambda_d\in\R_+$. Under \eqref{inegsqrt},
   $$\sum_{k=1}^d\lambda_k\xi_k^*\Sigma\xi_k\le \sum_{k=1}^d\lambda_k\left(\sum_{i=1}^n p_i\sqrt{\xi_k^*\Sigma_i\xi_k}\right)^2\le \sum_{i=1}^n p_i\sum_{k=1}^d\lambda_k\xi_k^*\Sigma_i\xi_k$$
   so that the integral of the convex quadratic function under ${\cal N}_d(0,\Sigma)$ is smaller than the convex combination with weights $p_i$ of the integrals under ${\cal N}_d(0,\Sigma_i)$ which, when $\sum_{i=1}^np_ix_i=0$, is equal to the convex combination with weights $p_i$ of the integrals under ${\cal N}_d(x_i,\Sigma_i)$ minus $\sum_{i=1}^n p_i\sum_{k=1}^d\lambda_k(\xi_k^*x_i)^2$.
 \end{Remark}
\begin{Remark}\label{remconsig}
  For given $\Sigma_1,\cdots,\Sigma_n\in{\cal S}_+(d)$ and $(p_1,\cdots,p_n)\in(0,1)^n$ such that $\sum_{i=1}^n p_i=1$ the set of matrices $\Sigma\in{\cal S}_+(d)$ which satisfy \eqref{inecov} and the set of matrices $\Sigma\in{\cal S}_+(d)$ which satisfy \eqref{inegsqrt} are convex. It is not clear whether the set of matrices $\Sigma$ that satisfy \eqref{compccgg} also is convex. The convexity would be clear  if the inequality ${\cal N}_d(0,p_1\Sigma_1+(1-p_1)\Sigma_2)\le_{cx} p_1{\cal N}_1(0,\Sigma_1)+(1-p_1){\cal N}_1(0,\Sigma_2)$ was valid for any $p_1\in[0,1]$ and any $\Sigma_1,\Sigma_2\in{\cal S}_+(d)$. Unfortunately, by \eqref{inegsqrt}, this inequality implies that $\forall \xi\in\R^d$, $p_1(1-p_1)(\sqrt{\xi^*\Sigma_1\xi}-\sqrt{\xi^*\Sigma_2\xi})^2\le 0$ so that either $p_1(1-p_1)=0$ or $\Sigma_1=\Sigma_2$.
\end{Remark}

 \begin{Remark}\label{condmultm}
  Note that when any of the two conditions \eqref{inecov} and \eqref{inegsqrt} holds, then it also holds with $(\Sigma,(\Sigma_i)_{i\in\{1,\cdots,n\}},d)$ replaced by  $(M\Sigma M^*,(M\Sigma_iM^*)_{i\in\{1,\cdots,n\}},q)$ for any $M\in\R^{q\times d}$. This is obvious for \eqref{inegsqrt}. Concerning \eqref{inecov}, this follows from the fact that for ${\mathbb M}\in\R^{nq\times nd}$ with blocks ${\mathbb M}_{(i-1)q+1:iq,(i-1)d+1:id}=M$ for $i\in\{1,\cdots,n\}$ and all other entries equal to $0$, $({\mathbb M}\Gamma{\mathbb M}^*)_{(i-1)q+1:iq,(i-1)q+1:iq}=M\Gamma_{(i-1)d+1:id,(i-1)d+1:id} M^*$  and $$(p_1I_q,\cdots,p_nI_q){\mathbb M}=(p_1M,\cdots,p_nM)=M(p_1I_d,\cdots,p_n I_d).$$ 
As a consequence, when $M\in{\cal I}(d)$,  any of the two conditions under consideration holds if and only if it holds with $(\Sigma,(\Sigma_i)_{i\in\{1,\cdots,n\}})$ replaced by $(M\Sigma M^*,(M\Sigma_iM^*)_{i\in\{1,\cdots,n\}})$.% Note that when any of the three conditions \eqref{inecov}, \eqref{compccgg} and \eqref{inegsqrt} holds, then it also holds with $(\Sigma,(\Sigma_i)_{i\in\{1,\cdots,n\}})$ replaced by  $(M\Sigma M^*,(M\Sigma_iM^*)_{i\in\{1,\cdots,n\}})$ for any $M\in\R^{q\times d}$. This is obvious for \eqref{inegsqrt} and this comes from the convexity of $\R^d\ni x\mapsto \varphi(Mx)$ for any convex function $\varphi:\R^q\to\R$ for \eqref{compccgg}. Concerning \eqref{inecov}, this follows from the fact that for ${\mathbb M}\in\R^{nq\times nd}$ with blocks ${\mathbb M}_{(i-1)q+1:id,(i-1)d+1:id}=M$ for $i\in\{1,\cdots,n\}$ and all other entries equal to $0$, $({\mathbb M}\Gamma{\mathbb M}^*)_{(i-1)q+1:id,(i-1)q+1:id}=M\Gamma_{(i-1)d+1:id,(i-1)d+1:id} M^*$  and $$(p_1I_d,\cdots,p_nI_q){\mathbb M}=(p_1M,\cdots,p_nM)=MA.$$ 
% As a consequence, when $M\in\R^{d\times d}$ is non singular,  any of the three conditions under consideration holds if and only if it holds with $(\Sigma,(\Sigma_i)_{i\in\{1,\cdots,n\}})$ replaced by $(M\Sigma M^*,(M\Sigma_iM^*)_{i\in\{1,\cdots,n\}})$.
\end{Remark}
\begin{Remark}\label{remcorrel}Let $\Lambda={\rm diag}(\lambda_1,\hdots,\lambda_d)$ where $\lambda_1,\cdots,\lambda_d\in\R\setminus\{0\}$. We remark that \eqref{correl} holds if and only if it holds with $M$ replaced by its product by $\Lambda$.  For the necessary condition, we decompose $\Lambda$ into the commutative product of $S={\rm diag}({\rm sign}(\lambda_1),\hdots,{\rm sign}(\lambda_d))$ and $|\Lambda|={\rm diag}(|\lambda_1|,\hdots,|\lambda_d|)$ and remark that $M\Sigma M^*\le DCD$ implies that $\Lambda M\Sigma M^*\Lambda\le \Lambda D C D\Lambda=(|\Lambda| D) SC S(|\Lambda| D)$ where $|\Lambda| D$ is diagonal with diagonal entries
  \begin{align*}
   (|\Lambda| D)_{kk}=|\lambda_k|D_{kk}=\sum_{i=1}^np_i\sqrt{\lambda_k^2(M\Sigma_i M^*)_{kk}}=\sum_{i=1}^np_i\sqrt{(\Lambda M\Sigma_i M^*\Lambda)_{kk}}\mbox{ for }k\in\{1,\cdots,d\}
  \end{align*}
  and, for $i\in\{1,\cdots,n\}$, $SCS$ is a correlation matrix associated with $\Lambda M\Sigma_i(\Lambda M)^*$ when $C$ is a correlation matrix associated with $M\Sigma_iM^*$. Since $\Lambda^{-1}={\rm diag}\left(\frac 1{\lambda_1},\hdots,\frac 1{\lambda_d}\right)$ where $\frac 1{\lambda_1},\cdots,\frac 1{\lambda_d}\in\R\setminus\{0\}$, the sufficient condition is proved in the same way.
 % Let us remark that \eqref{correl} holds if and only if it holds with $M$ replaced by its product by any diagonal matrix $\Lambda={\rm diag}(\lambda_1,\hdots,\lambda_d)$ with $\lambda_1,\cdots,\lambda_d>0$ (real?). The sufficient condition follows from the choice $\lambda_1=\cdots=\lambda_d=1$. For the necessary condition, we remark that $M\Sigma M^*\le DCD$ implies that $\Lambda M\Sigma M^*\Lambda\le \Lambda D C D\Lambda=(\Lambda D) C (\Lambda D)$ where $\Lambda D$ is diagonal with diagonal entries
 %  \begin{align*}
 %   (\Lambda D)_{kk}=\lambda_kD_{kk}=\sum_{i=1}^np_i\sqrt{\lambda_k^2(M\Sigma_i M^*)_{kk}}=\sum_{i=1}^np_i\sqrt{(\Lambda M\Sigma_i M^*\Lambda)_{kk}}\mbox{ for }k\in\{1,\cdots,d\}.
 %  \end{align*}
\end{Remark}

\subsection{Proofs of the main results}\label{sec:proof}
The proofs rely on the next result which is stated for instance in \cite[Remark 3.1]{JP}.
\begin{Lemma}\label{lemcnsord}
For $M,N\in{\cal S}_+(d)$, ${\cal N}_d(0,M)\le_{cx}{\cal N}_d(0,N)\Leftrightarrow M\le N$. 
\end{Lemma}

We give the very short proof of this lemma for the sake of completeness.

\noindent{\bf Proof.}
Let $M\le N$ and $X\sim{\cal N}_d(0,M)$ be independent of $Y\sim{\cal N}_d(0,N-M)$. Then $X+Y\sim{\cal N}_d(0,N)$ and since $\E(X+Y|X)=X$, the inequality ${\cal N}_d(0,M)\le_{cx}{\cal N}_d(0,N)$ follows from Jensen's inequality. Conversely, since for $Z\sim{\cal N}_d(0,\Sigma)$ and $\xi\in\R^d$, $\E((\xi^*Z)^2)=\xi^*\Sigma\xi$ where $\R^d\ni x\mapsto(\xi^*x)^{2}$ is convex, ${\cal N}_d(0,M)\le_{cx}{\cal N}_d(0,N)\Rightarrow \forall \xi\in\R^d,\;\xi^*M\xi\le \xi^*N\xi$.
\hfill$\Box$

Let us give a very short proof of Proposition \ref{propcnsconvcomble} (which also is a special case of the second assertion in Proposition \ref{propradno} applied with $q=d$, $Z\sim{\cal N}_d(0,I_d)$, $\sigma=\Sigma^{1/2}$ and $\sigma_i=\Sigma_i^{1/2}$ for $i\in\{1,\cdots,n\}$) before that of Theorem \ref{proprinc}.

\noindent{\bf Proof of Proposition \ref{propcnsconvcomble}.}
The sufficient condition follows from Lemma \ref{lemcnsord}. For the necessary condition, we suppose that $\sum_{i=1}^np_i{\cal N}_d(0,\Sigma_i)\le_{cx}{\cal N}_d(0,\Sigma)$. For the convex function $\varphi(x)=e^{\lambda \xi^*x}$ with $\lambda\in\R$ and $\xi\in\R^d$, we deduce that $\sum_{i=1}^np_ie^{\frac{\lambda^2}{2}\xi^*\Sigma_i\xi}\le e^{\frac{\lambda^2}{2}\xi^*\Sigma\xi}$ so that for each $i\in\{1,\cdots,n\}$, $\lambda^2(\xi^*\Sigma_i\xi-\xi^*\Sigma\xi)+\ln(p_i)\le 0$. By letting $|\lambda|\to\infty$, we deduce that $\xi^*\Sigma_i\xi\le \xi^*\Sigma\xi$.

% Since for $X\sim{\cal N}_d(0,M)$, $\xi\in\R^d$ and $k\in{\mathbb N}$, $\E(|\xi^*X|^{2k})=\frac{(\xi^*M\xi)^k(2k)!}{2^kk!}$ and $\R^d\ni x\mapsto|\xi^*x|^{2k}$ is convex, we deduce that $$\forall \xi\in\R^d,\;\left(\sum_{i=1}^np_i (\xi^*\Sigma_i\xi)^k\right)^{1/k}\le \xi^*\Sigma\xi.$$ Letting $k\to\infty$, we deduce that
  % $$\forall \xi\in\R^d,\;\max_{1\le i\le n}(\xi^*\Sigma_i\xi)\le\xi^*\Sigma\xi$$
  % so that $\Sigma_i\le \Sigma$ for all $i\in\{1,\cdots,n\}$.
  \hfill$\Box$

\noindent{\bf Proof of Theorem \ref{proprinc}.}
In this proof, we set $A=(p_1I_d,\cdots,p_n I_d)\in\R^{d\times nd}$.

Let us assume \eqref{correl} and let $B\in\R^{nd\times d}$ be defined by \begin{equation}
   B_{(i-1)d+1:id,1:d}={\rm diag}(\sqrt{(M\Sigma_i M^*)_{11}},\cdots,\sqrt{(M\Sigma_i M^*)_{dd}})\mbox{ for }i\in\{1,\cdots,n\}.\label{defB}
 \end{equation} Then $BCB^*\in{\cal S}_+(nd)$ and $(BCB^*)_{(i-1)d+1:id,(i-1)d+1:id}=M\Sigma_i M^*$. Moreover, for $A=(p_1I_d,\cdots,p_n I_d)\in\R^{d\times nd}$, $AB=D$ so that $ABCB^*A^*=(AB)C(AB)^*=DCD\ge M\Sigma M^*$. With Remark \ref{condmultm}, we conclude that \eqref{inecov} holds.

 Let us now assume \eqref{inecov}. By Lemma \ref{lemcnsord}, the inequality in this condition implies that
\begin{align}
   {\cal N}_d(0,\Sigma)\le_{cx}{\cal N}_d(0,A\Gamma A^*).\label{compgg}
\end{align}
Let $X\sim{\cal N}_{nd}(0,\Gamma)$. Then $AX=\sum_{i=1}^n p_iX_{(i-1)d+1:id}\sim{\cal N}_d(0,A\Gamma A^*)$. Let $x_1,\cdots,x_n\in\R^d$ be such that $\sum_{i=1}^n p_ix_i=0$. Since for $\varphi:\R^d\to\R$ convex, $$\varphi\left(\sum_{i=1}^n p_iX_{(i-1)d+1:id}\right)=\varphi\left(\sum_{i=1}^n p_i(X_{(i-1)d+1:id}+x_i)\right)\le \sum_{i=1}^n p_i\varphi(X_{(i-1)d+1:id}+x_i),$$ one has $${\cal N}_d(0,A\Gamma A^*)\le_{cx}\sum_{i=1}^n p_i{\cal N}_d(x_i,\Sigma_i),$$ inequality which, combined with \eqref{compgg}, implies \eqref{compccgg}.

For the convex function $\varphi(x)=\sqrt{\pi/2}|\xi^*x|$ where $\xi\in\R^d$, ${\cal N}_d(0,\Sigma)\le_{cx}\sum_{i=1}^np_i{\cal N}_d(0,\Sigma_i)$ implies \eqref{inegsqrt}.

In each of the four settings given at the end of the statement, let us prove that \eqref{inegsqrt} implies \eqref{correl} to check the equivalence of the four conditions.

When the matrices $\Sigma_i,\;i\in\{1,\cdots,n\}$ all are colinear i.e. there exist $\Theta\in{\cal S}_+(d)$ and $\sigma_1,\cdots,\sigma_d\in\R_+$ such that $\Sigma_i=\sigma_i\Theta$ for $i\in\{1,\cdots,n\}$, then $\sqrt{\xi^*\Sigma_i\xi}\sqrt{\xi^*\Sigma_j\xi}=\sqrt{\sigma_i\sigma_j}\times\xi^* \Theta\xi$ for all $i,j\in\{1,\cdots,n\}$ and $\xi\in\R^d$ so that $\left(\sum_{i=1}^np_i\sqrt{\xi^*\Sigma_i\xi}\right)^2=\left(\sum_{i=1}^np_i\sqrt{\sigma_i}\right)^2\xi^*\Theta\xi$ and \eqref{inegsqrt} implies $\Sigma\le \left(\sum_{i=1}^np_i\sqrt{\sigma_i}\right)^2\Theta$. Let $C$ be a correlation matrix associated with $\Theta$ and therefore  with each the matrices $\Sigma_i$, $i\in\{1,\cdots,n\}$. Then one has
$\Theta={\rm diag}\left(\sqrt{\Theta_{11}},\cdots,\sqrt{\Theta_{dd}}\right)C{\rm diag}\left(\sqrt{\Theta_{11}},\cdots,\sqrt{\Theta_{dd}}\right)$ and, since 
$(\Sigma_i)_{kk}=\sigma_i\Theta_{kk}$ for $i\in\{1,\cdots,n\}$ and $k\in\{1,\cdots,d\}$, $\left(\sum_{i=1}^np_i\sqrt{\sigma_i}\right)^2\Theta$ is equal to $${\rm diag}\left(\sum_{i=1}^np_i\sqrt{(\Sigma_i)_{11}},\cdots,\sum_{i=1}^np_i\sqrt{(\Sigma_i)_{dd}}\right)C{\rm diag}\left(\sum_{i=1}^np_i\sqrt{(\Sigma_i)_{11}},\cdots,\sum_{i=1}^np_i\sqrt{(\Sigma_i)_{dd}}\right),$$
so that \eqref{correl} holds with $M=I_d$.

When there exists $N\in{\cal I}(d)$ such that $N\Sigma_iN^*N\Sigma_jN^*=0$ for $i,j\in\{1,\cdots,n\}$ such that $i\ne j$, then the matrices $N\Sigma_iN^*\in{\cal S}_+(d),i\in\{1,\cdots,n\}$ commute and therefore are co-diagonalizable and there exist $O\in{\cal O}(d)$, $\lambda_1,\cdots\lambda_d\ge 0$ and $0\le k_1\le k_2\le \cdots\le k_n\le d$ such that $ON\Sigma_i(ON)^*={\rm diag}(0,\cdots,0,\lambda_{k_{i-1}+1},\cdots,\lambda_{k_i},0,\cdots,0)$ for $i\in\{1,\cdots,n\}$ (with conventions $k_0=0$ and that the matrix ${\rm diag}(0,\cdots,0,\lambda_{k_{i-1}+1},\cdots,\lambda_{k_i},0,\cdots,0)=0$ if $k_{i-1}=k_i$). Choosing $\xi=(ON)^*\zeta$ in \eqref{inegsqrt} with $\zeta\in\R^d$ with coordinates with indices outside $\{k_{i-1}+1,\cdots,k_i\}$ equal to $0$, we obtain that $(ON\Sigma (ON)^*)_{k_{i-1}+1:k_i,k_{i-1}+1:k_i}\le p_i^2{\rm diag}(\lambda_{k_{i-1}+1},\cdots,\lambda_{k_i})$. Hence $ON\Sigma(ON)^*\le\widetilde\Sigma$ where $\widetilde\Sigma\in{\cal S}_+(d)$ denotes the matrix with diagonal blocks $\widetilde\Sigma_{k_{i-1}+1:k_i,k_{i-1}+1:k_i}=p_i^2{\rm diag}(\lambda_{k_{i-1}+1},\cdots,\lambda_{k_i})$ and all other entries equal to those of $ON\Sigma (ON)^*$. Let $C$ be a correlation matrix associated with $\widetilde\Sigma$. Then for $i\in\{1,\cdots,n\}$, since $ON\Sigma_i(ON)^*$ coincides with $\widetilde\Sigma$ on the diagonal block with indices between $k_{i-1}+1$ and $k_i$ and has all other entries equal to $0$, $C$ is associated with $ON\Sigma_i (ON)^*$. Moreover, for $k\in\{1,\cdots,d\}$,
$$\sqrt{\widetilde\Sigma_{kk}}=\sqrt{\sum_{i=1}^n 1_{\{k_{i-1}+1\le k\le k_i\}}p_i^2\lambda_k}=\sum_{i=1}^n 1_{\{k_{i-1}+1\le k\le k_i\}}p_i\sqrt{\lambda_k}=\sum_{i=1}^np_i\sqrt{(ON\Sigma_i (ON)^*)_{kk}},$$
so that \eqref{correl} holds with $M=ON$.

Let us finally suppose the existence of $(M,C)\in{\cal I}(d)\times{\cal C}(d)$ such that $C$ is associated with all the matrices $\widehat\Sigma,M\Sigma_1 M^*,\cdots, M\Sigma_nM^*$. This condition is satisfied when there exists $N\in{\cal I}(d)$ such that the matrices $N\Sigma N^*,N\Sigma_1 N^*,\cdots,N\Sigma_n N^*$ commute and therefore are co-diagonalizable. Indeed, for $M=ON$ where $O\in{\cal O}(d)$ is such that $ON\Sigma (ON)^*,ON\Sigma_1 (ON)^*,\cdots,ON\Sigma_n(ON)^*$ are diagonal, the identity correlation matrix $C=I_d$ is associated with all the matrices $\widehat \Sigma,M\Sigma_1 M^*,\cdots,M\Sigma_nM^*$.
 Choosing $\xi$ as the $k$-th row of $M$ in \eqref{inegsqrt}, we obtain that $(M\Sigma M^*)_{kk}\le\left(\sum_{i=1}^n p_i\sqrt{(M\Sigma_i M^*)_{kk}}\right)^2$ for $k\in\{1,\cdots,d\}$. Hence $M\Sigma M^*\le \widehat \Sigma=DCD$ 
so that \eqref{correl} holds.\hfill$\Box$
 \begin{Remark}
   According to Strassen's theorem, the convex ordering between two probability distributions is equivalent to the existence of a martingale coupling between the smaller and the larger. Under \eqref{inecov}, when $\sum_{i=1}^n p_ix_i=0$, such a coupling between ${\cal N}_d(0,\Sigma)$ and $\sum_{i=1}^np_i{\cal N}_d(x_i,\Sigma_i)$ is given ${\cal N}_d(0,\Sigma)(dx)P(x,dy)$ with $P(x,dy)$ the martingale Markov kernel defined by $$P(x,dy)=\int_{(z,w)\in\R^d\times\R^{nd}}\sum_{i=1}^n p_i\delta_{w_{(i-1)d+1:id}+x_i}(dy)Q(z,dw){\cal N}_d(x,A\Gamma A^*-\Sigma)(dz)$$
   where $Q(z,dw)$ denotes the conditional law of $X\sim{\cal N}_{nd}(0,\Gamma)$ given $\sum_{i=1}^n p_iX_{(i-1)d+1:id}=z$ and $A=(p_1I_d,\cdots,p_nI_d)$. Note that $P(x,dy)$ is a Gaussian mixture.
  % \textcolor{red}{calculer ses caractéristiques n'est pas facile : d'abord ${\cal L}(X_{(1:(i-1)d)(id+1:nd)}|X_{(i-1)d+1:id)})={\cal N}_{(n-1)d}(0,\Theta^{(i)})$ où $\Theta^{(i)}=\Gamma_{(1:(i-1)d)(id+1:nd),(1:(i-1)d)(id+1:nd)}-\Gamma_{(1:(i-1)d)(id+1:nd),(i-1)d+1:id}\Sigma_i^{-1}\Gamma_{(1:(i-1)d)(id+1:nd),(i-1)d+1:id})$ puis ${\cal L}(AX|X_{(i-1)d+1:id)})={\cal N}_{d}(p_iX_{(i-1)d+1:id)},A_{-i}\Theta^{(i)}A_{-i}^*$ où $A_{-i}=(p_1I_d,\cdots,p_{i-1}I_d,p_{i+1}I_d,\cdots,p_nI_d)$ puis ${\cal L}(X_{(i-1)d+1:id)}|AX)$ puis ${\cal L}(X_{(i-1)d+1:id)}|x)$ les deux derniers en utilisant le thm MAP534}
\end{Remark}
There is of course no need to go through convex ordering of Gaussian mixtures to check that \eqref{inecov} implies \eqref{inegsqrt}, as confirmed by the next Lemma. 

\begin{Lemma}\label{lemcovfsq}
  Let $n\ge 2$, $\Sigma,\Sigma_1,\cdots,\Sigma_n\in{\cal S}_+(d)$ and $(p_1,\cdots,p_n)\in(0,1)^n$ such that $\sum_{i=1}^n p_i=1$.
For $\Gamma\in\R^{nd\times nd}$ and $i,j\in\{1,\cdots,n\}$, we denote $\Gamma_{(ij)}=(\Gamma_{(i-1)d+k,(j-1)d+\ell})_{1\le k,\ell\le d}\in\R^{d\times d}$. Then \eqref{inecov} implies
\begin{align}
 \exists\Gamma\in\R^{nd\times nd}&\mbox{symmetric such that }\Sigma\le (p_1I_d,\cdots,p_n I_d)\Gamma (p_1I_d,\cdots,p_n I_d)^*,\;\Gamma_{(ii)}=\Sigma_i\mbox{ for }1\le i\le n,\;\notag\\&
  \mbox{ and }\left(\begin{array}{cc}\Gamma_{(ii)} & \Gamma_{(ij)}
   \\\Gamma_{(ji)} &\Gamma_{(jj)}
\end{array}\right)\in{\cal S}_+(2d)\mbox{ for }1\le i<j\le n\label{inecovf},
\end{align}which in turn implies  \eqref{inegsqrt}. 
\end{Lemma}

\noindent{\bf Proof.}
Since $\Gamma\in{\cal S}_+(nd)$ implies that $\left(\begin{array}{cc}\Gamma_{(ii)} & \Gamma_{(ij)}
   \\\Gamma_{(ji)} &\Gamma_{(jj)}
                                                        \end{array}\right)\in{\cal S}_+(2d)$ for $1\le i<j\le n$, one has \eqref{inecov}$\Rightarrow$\eqref{inecovf}.
Under \eqref{inecovf}, for $\xi\in\R^d$ and $1\le i<j\le n$, one has $\Gamma_{(ji)}=\Gamma_{(ij)}^*$ and therefore $\xi^*\Gamma_{(ji)}\xi=\xi^*\Gamma_{(ij)}\xi$ by symmetry of $\Gamma$ and
$$\left(\begin{array}{cc} \xi^*\Sigma_i\xi & \xi^*\Gamma_{(ij)}\xi
   \\ \xi^*\Gamma_{(ij)}\xi & \xi^*\Sigma_j\xi\end{array}\right)=\left(\begin{array}{cc} \xi& 0 \\ 0 &\xi\end{array}\right)^* \left(\begin{array}{cc}\Gamma_{(ii)} & \Gamma_{(ij)}
   \\\Gamma_{(ji)} &\Gamma_{(jj)}
                                                        \end{array}\right)\left(\begin{array}{cc} \xi & 0\\0& \xi\end{array}\right)\in{\cal S}_+(2)$$
so that $|\xi^*\Gamma_{(ij)}\xi  |\le\sqrt{\xi^*\Sigma_i\xi\times\xi^*\Sigma_j\xi}$. As a consequence, 
 \begin{align*}
   \xi^*\Sigma\xi&\le \xi^*(p_1I_d,\cdots,p_nI_d)\Gamma (p_1I_d,\cdots,p_nI_d)^*\xi=\sum_{i,j=1}^n p_ip_j\xi^*\Gamma_{(ij)}\xi\\&\le \sum_{i,j=1}^n p_ip_j\sqrt{\xi^*\Sigma_i\xi\times\xi^*\Sigma_j\xi}=\left(\sum_{i=1}^n p_i\sqrt{\xi^*\Sigma_i\xi}\right)^2.\end{align*} 
\hfill$\Box$
\begin{Remark}We deduce that under equivalence of \eqref{inecov} and \eqref{inegsqrt} and in particular in the four settings given at the end of Theorem \ref{proprinc},  then \eqref{inecovf} and \eqref{compccgg} also are equivalent.

  Note that $\left(\begin{array}{cc}\Gamma_{(ii)} & \Gamma_{(ij)}
   \\\Gamma_{(ji)} &\Gamma_{(jj)}
                                                        \end{array}\right)\in{\cal S}_+(2d)$ for $1\le i<j\le n$ does not imply that $\Gamma\in{\cal S}_+(nd)$ when $n\ge 3$. As a consequence, it is very likely that \eqref{inecov} is stronger than \eqref{inegsqrt} when $n\ge 3$, even if we have not been able to provide an example where \eqref{inegsqrt} holds but \eqref{inecov} does not in reason of the great flexibility offered by the choice of $\Gamma$ in this condition (once the symmetry and the equality of the diagonal blocks of $\Gamma$ to the matrices $\Sigma_i$ are taken into account, there are still $\frac{n(n-1)}{2}d^2$ entries to choose). \end{Remark}
                                                    \subsection{Examples}
 According to the next example, \eqref{inecov} does not imply \eqref{correl}.    \begin{Example}
   Let $n=d=2$, $p_1=p_2=\frac 12$, $\Sigma_1=\left(\begin{array}{cc} 4 & 0
   \\ 0 & 4
                \end{array}\right)$, $\Sigma_2=\left(\begin{array}{cc} 4 & 0
   \\ 0 & 4\lambda^2
                \end{array}\right)$ and $\Sigma=\left(\begin{array}{cc} 2 & 1+\lambda
   \\ 1+\lambda & 1+\lambda^2
                                                      \end{array}\right)$ where $\lambda\in\R\setminus\{-1,1\}$ so that $\Sigma_1\ne \Sigma_2$.
% For $\xi=\left(\begin{array}{c} \xi_1
%    \\ \xi_2
%                 \end{array}\right)\in\R^2$, \begin{align*}
%                   \left(\frac{\sqrt{\xi^*\Sigma_1\xi}+\sqrt{\xi^*\Sigma_2\xi}}2\right)^2&=\left(\sqrt{\xi_1^2+\xi_2^2}+\sqrt{\xi_1^2+\lambda^2\xi_2^2}\right)^2
%                                                                                           =2\xi_1^2+(1+\lambda^2)\xi_2^2+2\sqrt{(1+\lambda)^2\xi_1^2\xi_2^2+(\xi^2_1-\lambda\xi^2_2)^2}\\&\ge  2\xi_1^2+(1+\lambda^2)\xi_2^2+2(1+\lambda)|\xi_1\xi_2|\ge \xi^*\Sigma\xi \end{align*} so that \eqref{inegsqrt} is satisfied. 
                                                                                        Condition \eqref{inecov} is satisfied with $\Gamma=4\left(\begin{array}{cccc} 1 & 0 &0 &\lambda
                                                                                                                                                                                                                                                                                                                                                                              \\ 0 & 1 &1&0\\0 & 1 &1&0\\\lambda & 0 &0 &\lambda^2
                \end{array}\right)$.

The identity correlation matrix $I_2$ is associated with both matrices $\Sigma_1$ and $\Sigma_2$. The matrix $D$ defined in \eqref{correl} with $M=I_2$ writes $D={\rm diag}\left(\frac{\sqrt{(\Sigma_1)_{11}}+\sqrt{(\Sigma_2)_{11}}}2,\frac{\sqrt{(\Sigma_1)_{22}}+\sqrt{(\Sigma_2)_{22}}}2\right)= {\rm diag}(2,1+\lambda)$.                                            
One has $DI_2D=\left(\begin{array}{cc} 4 & 0
   \\ 0 & (1+\lambda)^2
                \end{array}\right)$ and $DI_2D-\Sigma=\left(\begin{array}{cc} 2 & -(1+\lambda)
   \\ -(1+\lambda) & 2\lambda
                                                            \end{array}\right)$ with determinant ${\rm det}(DI_2D-\Sigma)=-(1-\lambda)^2$ so that $\Sigma$ is not smaller than $DI_2D$.
In order to assess whether the more general condition \eqref{correl} is satisfied, let us now find the matrices $M=\left(\begin{array}{cc} a & b
                                                     \\ c & d\end{array}\right)\in\R^{2\times 2}$ non-singular i.e. with entries satisfying $ad-bc\neq 0$ which are such that the same correlation matrix $C$ is associated with both $M\Sigma_1M^*=4\left(\begin{array}{cc} a^2+b^2 & ac+bd
                                                     \\ ac+bd & c^2+d^2\end{array}\right)$ and $M\Sigma_2M^*=4\left(\begin{array}{cc} a^2+\lambda^2b^2 & ac+\lambda^2bd
                                                                                                                      \\ ac+\lambda^2bd & c^2+\lambda^2d^2\end{array}\right)$ or equivalently
                                                                                                                  $$\frac{ac+bd}{\sqrt{(a^2+b^2)(c^2+d^2)}}=\frac{ac+\lambda^2bd}{\sqrt{(a^2+\lambda^2b^2)(c^2+\lambda^2d^2)}}.$$
                                                                                                                  Multiplying both sides by $\sqrt{(a^2+b^2)(c^2+d^2)(a^2+\lambda^2b^2)(c^2+\lambda^2d^2)}$ and then squaring we deduce that
\begin{align*}
   &(ac+bd)^2((ac+\lambda^2bd)^2+\lambda^2(ad-bc)^2)=(ac+\lambda^2bd)^2((ac+bd)^2+(ad-bc)^2)\\\Leftrightarrow &(ad-bc)^2(\lambda^2-1)((ac)^2-(\lambda bd)^2)=0\Leftrightarrow |ac|=|\lambda bd|.
\end{align*}
Hence a necessary condition is $|ac|=|\lambda bd|$. Conversely, when $ac=|\lambda| bd$, the correlation matrix % $C=\left(\begin{array}{cc} 1 & \frac{(1+\lambda)bd}{\sqrt{(a^2+b^2)(c^2+d^2)}} \\\frac{(1+\lambda)bd}{\sqrt{(a^2+b^2)(c^2+d^2)}} & 1\end{array}\right)$
$C\in\R^{2\times 2}$ with extra-diagonal entries equal to $\frac{(1+|\lambda|)bd}{\sqrt{(a^2+b^2)(c^2+d^2)}}$  is associated with both $M\Sigma_1M^*$ and $M\Sigma_2M^*$ while when $ac=-|\lambda| bd$ the respective extra-diagonal entries $\frac{(1-|\lambda|)bd}{\sqrt{(a^2+b^2)(c^2+d^2)}}$ and  $\frac{(|\lambda|-1)bd}{\sqrt{(a^2+b^2)(c^2+d^2)}}$ of the correlation matrices associated with $M\Sigma_1M^*$ and $M\Sigma_2M^*$ differ unless $bd=0$. As a consequence, the same correlation matrix is associated with both $M\Sigma_1M^*$ and $M\Sigma_2M^*$ if and only if $ac=|\lambda| bd$, a condition which, together with $ad-bc\ne 0$, implies that when moreover $bd=0$ then either $a=d=0$ or $b=c=0$. 
Let us now suppose that $\lambda=0$. The necessary and sufficient condition writes $ac=0$. When $c=0$ (resp. $a=0$), then $ad\ne 0$ (resp. $bc\ne 0$) and according to Remark \ref{remcorrel} applied with $\Lambda={\rm diag}(\frac 1 a,\frac 1 d)$ (resp. $\Lambda={\rm diag}(\frac 1 b,\frac 1 c)$), \eqref{correl} holds if and only if it holds with $M=\left(\begin{array}{cc} 1 & x \\0& 1\end{array}\right)$ (resp. $M=\left(\begin{array}{cc} 0& 1\\1 & x \end{array}\right)$) where $x=\frac b a$ (resp. $x=\frac d c$). Since \eqref{correl} is clearly preserved under exchange of the rows of $M$, it is enough to deal with $M=\left(\begin{array}{cc} 1 & x \\0& 1\end{array}\right)$. Then the correlation matrix associated to both $M\Sigma_1M^*=\left(\begin{array}{cc} 4(1+x^2) & 4x \\4x & 4\end{array}\right)$ and $M\Sigma_2M^*=\left(\begin{array}{cc} 4 & 0 \\0 & 0\end{array}\right)$ is $C=\left(\begin{array}{cc} 1 & \frac{x}{\sqrt{1+x^2}} \\\frac{x}{\sqrt{1+x^2}}& 1\end{array}\right)$. On the other hand, $M\Sigma M^*=\left(\begin{array}{cc} 2+2x+x^2 & 1+x \\1+x & 1\end{array}\right)$ while the matrix $D$ defined in \eqref{correl} writes $D={\rm diag}(1+\sqrt{1+x^2},1)$ so that $DCD=\left(\begin{array}{cc} (1+\sqrt{1+x^2})^2 & x+\frac{x}{\sqrt{1+x^2}} \\x+\frac{x}{\sqrt{1+x^2}}& 1\end{array}\right)$ and $DCD-M\Sigma M^*=\left(\begin{array}{cc} 2(\sqrt{1+x^2}-x) & \frac{x}{\sqrt{1+x^2}}-1 \\\frac{x}{\sqrt{1+x^2}}-1 & 0\end{array}\right)$ is not positive semi-definite. Hence \eqref{correl} does not hold while \eqref{inecov} holds. % As a consequence, there exists some $\varepsilon>0$ such that this remains true for $\lambda\in(-\varepsilon,\varepsilon)$. Otherwise, one could find a sequence $(\lambda_n)_{n\in\N}$ converging to $0$ as $n\to\infty$ such that \eqref{correl} holds for $\lambda =\lambda_n$ with $M=M_n$. Since \eqref{correl} is clearly preserved by division of $M$ by the maximum of its absolute entries, up to extracting a subsequence, one can suppose that $M_n$ converges to some limit $M_\infty$.
\end{Example}
In the next example with $n=d=2$, we are able to exhibit the covariance matrices $\Sigma$ with equal diagonal coefficients such that \eqref{inegsqrt} holds. We also check that they satisfy \eqref{inecov} and therefore \eqref{compccgg}.
\begin{Example}\label{exdiag2d}For $\Sigma_1=\left(\begin{array}{cc}8 & 0
   \\ 0 & 4
\end{array}\right)$, $\Sigma_2=\left(\begin{array}{cc}4 & 0
   \\ 0 & 8
                                     \end{array}\right)$ and $p_1=p_2=\frac 12$, let us check  that the set of covariance matrices $\Sigma\in{\cal S}_+(2)$ with equal diagonal coefficients such that ${\cal N}_2(0,\Sigma)\le_{cx}\frac{1}{2}{\cal N}_2(0,\Sigma_1)+\frac{1}{2}{\cal N}_2(0,\Sigma_2)$ is the set of matrices $\left(\begin{array}{cc} a & b\\b & a\end{array}\right)$ with $a\in[0,3]$ and $|b|\le a$ or $a\in[3,\frac{17}3]$ and $|b|\le 6-a$ or $a\in[\frac{17}3,3+2\sqrt{2}]$ and $|b|\le \sqrt{1-\frac{(a-3)^2}8}$. We first check that this is also the set of covariance matrices $\Sigma\in{\cal S}_+(2)$ with equal diagonal coefficients such that \eqref{inegsqrt} holds. We then check that this set is included in the set of covariance matrices $\Sigma\in{\cal S}_+(2)$ such that \eqref{inecov} holds and conclude by Theorem \ref{proprinc}.

                                  Any element of ${\cal S}_+(2)$ with equal diagonal coefficients writes $\Sigma=\left(\begin{array}{cc} a & b\\b & a\end{array}\right)$ with $a\ge 0$ and $b^2\le a^2$. One has $\frac{1+\alpha}4\Sigma_1+\frac {1+1/\alpha}{4}\Sigma_2=\left(\begin{array}{cc}3+2\alpha+\frac 1{\alpha} & 0
   \\ 0 & 3+\alpha+\frac 2{\alpha}
                                                                                                                                                                                                                                                                               \end{array}\right)$ and $\inf_{\alpha>0}\left(2\alpha+\frac 1{\alpha}\right) =2\sqrt{2}$. Hence \eqref{inegsqrtbis} and \eqref{inegsqrt} hold if and only if $a\le 3+2\sqrt{2}$ and \begin{equation}\forall \alpha>0,\;b^2\le \left(3-a+2\alpha+\frac 1{\alpha}\right)\left(3-a+\alpha+\frac 2{\alpha}\right)=(3-a)^2+1+g\left(\alpha+\frac 1\alpha\right),\label{bal1al}\end{equation}                                                                                                                                    where $g(x)=2x^2+3(3-a)x$ is non-decreasing on $[\frac{3(a-3)}4,+\infty)$. This is the place where the choice of equal diagonal coefficients plays a crucial role. Indeed, with different coefficients, the product in the right-hand side of \eqref{bal1al} would not depend on $(\alpha,\frac 1\alpha)$ only through the sum $\alpha+\frac 1\alpha$. The simple dependence through the sum permits explicit minimization. The function $(0,+\infty)\ni\alpha\mapsto \alpha+\frac{1}{\alpha}$ is minimal and equal to $2$ for $\alpha=1$. When $a\in[0,\frac{17}3]$, $\frac{3(a-3)}4\le 2$ and \eqref{inegsqrt} holds if and only if $b^2\le (3-a)^2+1+g(2)=(6-a)^2$, a constraint which is implied by $b^2\le a^2$ when $a\in[0,3]$. When $a\in[\frac{17}3,3+2\sqrt{2}]$, $\frac{3(a-3)}4\ge 2$ and \eqref{inegsqrt} holds if and only if $b^2\le (3-a)^2+1+g\left(\frac{3(a-3)}4\right)=1-\frac{(a-3)^2}8$.

For $a\in[\frac{17}3,3+2\sqrt{2}]$, the matrix $\Gamma=\left(\begin{array}{cc}\Sigma_1 & \Theta
                                                               \\ \Theta^* & \Sigma_2\end{array}\right)$ with $\Theta=\left(\begin{array}{cc}2(a-3) & \pm 8\sqrt{1-\frac{(a-3)^2}8}                                                                                        \\ \mp 4\sqrt{1-\frac{(a-3)^2}8} & 2(a-3)\end{array}\right)$ is positive semi-definite  since $$\Gamma=\left(\begin{array}{c} S_1 \\ S_2
                 \end{array}\right)\left(\begin{array}{c} S_1 \\ S_2
                 \end{array}\right)^*\mbox{ for }S_1=\left(\begin{array}{cc} {2\sqrt{2}} &0 \\ 0 &2
                                                                                                                                                                                                                                                                                                    \end{array}\right)\mbox{ and }S_2=\left(\begin{array}{cc}\frac{a-3}{\sqrt{2}} &    \mp 2\sqrt{1-\frac{(a-3)^2}8}                                                                                  \\ \pm 2\sqrt{2}\sqrt{1-\frac{(a-3)^2}8}    & a-3\end{array}\right).$$
                                                                         We deduce that for $a\in[\frac{17}{3},3+2\sqrt{2}]$, \eqref{inecov} is satisfied when  $b\in\left\{-\sqrt{1-\frac{(a-3)^2}8},\sqrt{1-\frac{(a-3)^2}8}\right\}$ and therefore when $b\in\left[-\sqrt{1-\frac{(a-3)^2}8},\sqrt{1-\frac{(a-3)^2}8}\right]$ by convexity (see Remark \ref{remconsig}). For $a\in[3,\frac{17}{3}]$, since $$\left(\begin{array}{cc}a & \pm(6-a)
   \\ \pm(6-a) & a\end{array}\right)\le \left(\begin{array}{cc}\frac{17}{3} & \pm(6-\frac{17}{3})
   \\ \pm(6-\frac{17}{3}) & \frac{17}{3} \end{array}\right)=\left(\begin{array}{cc}\frac{17}{3} & \pm\sqrt{1-\frac{(\frac{17}{3}-3)^2}8}
   \\ \pm\sqrt{1-\frac{(\frac{17}{3}-3)^2}8} & \frac{17}{3} \end{array}\right),$$ we deduce that \eqref{inecov} is satisfied when $b\in\{a-6,6-a\}$ and therefore when $|b|\le 6-a$ by convexity. For $a\in[0,3]$, since $\left(\begin{array}{cc}a & \pm a                                                                                       \\  \pm a   & a\end{array}\right)\le  \left(\begin{array}{cc} 3 & \pm 3                                                                                       \\  \pm 3  & 3\end{array}\right)=\left(\begin{array}{cc} 3 & \pm (6-3)                                                                                      \\  \pm (6-3)  & 3\end{array}\right)$, we deduce that \eqref{inecov} is satisfied when $b\in\{-a,a\}$ and therefore when $b\in[-a,a]$ by convexity.
\end{Example}

In the next example with $n=3$ and $p_1=p_2=p_3=\frac 13$, for matrices $\Sigma_1,\Sigma_2,\Sigma_3\in{\cal S}_+(2)$ well-chosen in view of the previous examples, we exhibit a set of matrices $\Sigma\in{\cal S}_+(2)$ such that \eqref{inecovf} and therefore \eqref{inegsqrt} (in view of Lemma \ref{lemcovfsq}) hold. This set is parametrized by two parameters. The matrix $\Gamma$ for which we check \eqref{inecovf} belongs to ${\cal S}_+(6)$ when the couple of parameters belongs to a curve.
\begin{Example}Let $n=3$, $d=2$, $\Sigma_1=\left(\begin{array}{cc}18 & 0\\0 & 9\end{array}\right)$, $\Sigma_2=\left(\begin{array}{cc}9 & 0\\0 & 9\end{array}\right)$, $\Sigma_3=\left(\begin{array}{cc}9 & 0\\0 & 18\end{array}\right)$ and $p_1=p_2=p_3=\frac 13$.
For $\Sigma\in{\cal S}_+(2)$ diagonal, \eqref{correl}, \eqref{inecov}, \eqref{inecovf}, \eqref{compccgg} and \eqref{inegsqrt} are equivalent according to Theorem \ref{proprinc} since $\Sigma,\Sigma_1,\Sigma_2$ and $\Sigma_3$ commute. The set of diagonal covariance matrices which satisfy all these conditions is $\{\Sigma\in{\cal S}_+(2)\mbox{ diagonal and such that }\Sigma_{11}\vee\Sigma_{22}\le 6+4\sqrt{2}\}$ since, by Remark \ref{remskew}, \eqref{inegsqrt} holds if and only if\begin{align*}
 \forall\alpha,\beta,\gamma >0,\;  \Sigma&\le\frac 1 9\left((1+1/\alpha +\gamma)\Sigma_1+(1+\alpha+\beta)\Sigma_2+(1+1/\beta+1/\gamma)\Sigma_3\right)\\&=\left(\begin{array}{cc}4+\alpha+\frac 2\alpha+\beta+\frac 1\beta+2\gamma+\frac 1\gamma & 0\\0 & 4+\alpha+\frac 1\alpha+\beta+\frac 2\beta+\gamma+\frac 2\gamma\end{array}\right)
  \end{align*}
  and  $\inf_{\alpha>0}\left(\alpha+\frac 2\alpha\right)=2\sqrt{2}$  while  $\inf_{\alpha>0}\left(\alpha+\frac 1\alpha\right)=2$. According to Example \ref{exdiag2d}, for $a\in[\frac{17}{3},3+2\sqrt{2}]$ and $\Theta(a)=\left(\begin{array}{cc}\frac{9}{2}(a-3) & \pm18\sqrt{1-\frac{(a-3)^2}8}                                                                                        \\ \mp 9\sqrt{1-\frac{(a-3)^2}8} & \frac{9}{2}(a-3)\end{array}\right)$,  the matrix $\left(\begin{array}{cc} \Sigma_1 & \Theta(a) \\ \Theta(a)^* &\Sigma_3\end{array}\right)$ belongs to ${\cal S}_+(4)$. For $x\in[\frac{2^{5/4}}{1+\sqrt{2}},1]$ ($\frac{2^{5/4}}{1+\sqrt{2}}\simeq 0.9852$) and $\tilde \Theta(x)=\left(\begin{array}{cc}9\sqrt{2}x & \pm 9\sqrt{2}\sqrt{1-x^2}\\\mp 9\sqrt{1-x^2} & 9x\end{array}\right)$, the matrix $\left(\begin{array}{cc}\Sigma_1 & \tilde \Theta(x)\\\tilde \Theta(x)^* & \Sigma_2\end{array}\right)$ belongs to ${\cal S}_+(4)$ since it writes $\left(\begin{array}{c} S_1 \\ S_2(x)
              \end{array}\right)\left(\begin{array}{c cc} S_1^* &S_2(x)^* 
                                      \end{array}\right)$ with $S_1=\left(\begin{array}{cc} 3\sqrt{2} & 0\\ 0 & 3\end{array}\right)$ and $S_2(x)=\left(\begin{array}{cc} 3x & \mp 3\sqrt{1-x^2}\\ \pm 3\sqrt{1-x^2} & 3x\end{array}\right)$. Note that this can be used like in Example \ref{exdiag2d} to check that for $x\in[\frac{2^{5/4}}{1+\sqrt{2}},1]$, $${\cal N}_2\left(0,\left(\begin{array}{cc} \frac{9}{4}(3+2\sqrt{2}x) &  b\\b & \frac{9}{2}(1+x)\end{array}\right)\right)\le_{cx} \frac 12{\cal N}_2(0,\Sigma_1)+\frac 12{\cal N}_2(0,\Sigma_2)\Longleftrightarrow|b|\le \frac 9 4(\sqrt{2}-1)\sqrt{1-x^2}.$$ The diagonal structure ensures that, in the analogous of the inequality
 \eqref{bal1al}, the right-hand side still only depends on $\alpha$ and $\frac 1\alpha$ through their sum $\alpha+\frac 1{\alpha}$, which permits explicit minimization.

 Hence $$\Sigma(a,x)=\left(\begin{array}{cc}1+a+2(1+\sqrt{2})x &\pm\left(\sqrt{1-\frac{(a-3)^2}8}+2(\sqrt{2}-1)\sqrt{1-x^2}\right)\\  \pm\left(\sqrt{1-\frac{(a-3)^2}8}+2(\sqrt{2}-1)\sqrt{1-x^2}\right) & 1+a+2(1+\sqrt{2})x\end{array}\right)$$ satisfies \eqref{inecovf} with $\Gamma$ equal to  $\Gamma(a,x)=\left(\begin{array}{ccc}\Sigma_1 & \tilde\Theta(x) & \Theta(a) \\\tilde\Theta(x)^* &\Sigma_2 &\hat\Theta(x)^* \\\Theta(a)^* &\hat\Theta(x) & \Sigma_3\end{array}\right)$ for $\hat\Theta(x)=\left(\begin{array}{cc}9x &\mp 9\sqrt{1-x^2}\\  \pm 9\sqrt{2}\sqrt{1-x^2} & 9\sqrt{2}x\end{array}\right)$. By Lemma \ref{lemcovfsq}, we deduce that for $(a,x)\in[\frac{17}{3},3+2\sqrt{2}]\times[\frac{2^{5/4}}{1+\sqrt{2}},1]$, $\Sigma(a,x)$ satisfies \eqref{inegsqrt}.  In general the matrix $\Gamma(a,x)$ does not belong to ${\cal S}_+(6)$. For instance, the smallest eigenvalue of the matrix $\Gamma(\frac{17}3,1)$ is close to $-2.58$. Nevertheless, for $a\in[\frac{17}{3},3+2\sqrt{2}]$, $f(a)=\sqrt{\frac{a+2\sqrt{2}-3}{4\sqrt{2}}}$ and $S_3(a)=\left(\begin{array}{cc}\frac{3}{2\sqrt{2}}(a-3) & \mp 3\sqrt{1-\frac{(a-3)^2}{8}}\\\pm 3\sqrt{2}\sqrt{1-\frac{(a-3)^2}{8}} & \frac{3}{2}(a-3)\end{array}\right)$, the matrix $\left(\begin{array}{c} S_1 \\ S_2(f(a))
              \\ S_3(a)   \end{array}\right)\left(\begin{array}{c cc} S_1^* &S_2(f(a))^* & S_3(a)^*
                                                  \end{array}\right)$ is equal to $\Gamma(a,f(a))$ (note that $4f(a)^2(1-f(a)^2)=1-\frac{(a-3)²}{8}$) which thus belongs to ${\cal S}_+(6)$. As a consequence, for $a\in[\frac{17}{3},3+2\sqrt{2}]$ the matrix $\Sigma=\left(\begin{array}{cc}\Sigma_{11}(a,f(a)) & b\\b & \Sigma_{11}(a,f(a))\end{array}\right)$ satisfies \eqref{inecov} when % $b=\pm\left(\sqrt{1-\frac{(a-3)^2}8}+2(\sqrt{2}-1)\sqrt{1-f(a)^2}\right)$ and therefore when 
                                                $|b|\le |\Sigma_{12}(a,f(a))|$ (the case with strict inequality is deduced from the case with equality by Remark \ref{remconsig}). Conditions \eqref{inecov} and therefore \eqref{compccgg} are also satisfied by any matrix in ${\cal S}_+(2)$ smaller than such a matrix. Let us now discuss the values $(a,x)\in[\frac{17}{3},3+2\sqrt{2}]\times[\frac{2^{5/4}}{1+\sqrt{2}},1]$ such that $\Sigma(a,x)$ is smaller than such a matrix, a property referred to "dominated" in what follows. The function $f$  is increasing on $[\frac{17}{3},3+2\sqrt{2}]$ from $f(\frac{17}3)=\frac{1+\sqrt{2}}{\sqrt{6}}\simeq 0.9856$ to $f(3+2\sqrt{2})=1$ and $a\mapsto\Sigma_{11}(a,f(a))$ also is increasing on $[\frac{17}{3},3+2\sqrt{2}]$. The function $[\frac{17}3,3+2\sqrt{2}]\ni a\mapsto a+\sqrt{1-\frac{(a-3)^2}8}$  equal to the term depending on $a$ in $\Sigma_{11}(a,x)+|\Sigma_{12}(a,x)|$ is decreasing since its derivative $1-\frac{a-3}{8\sqrt{1-\frac{(a-3)^2}{8}}}$ is equal to $0$ for $a=\frac{17}3$ and decreasing on $[\frac{17}3,3+2\sqrt{2})$. The function $[\frac{2^{5/4}}{1+\sqrt{2}},1]\ni x\mapsto (1+\sqrt{2})x+(\sqrt{2}-1)\sqrt{1-x^2}$ % =(1+\sqrt{2})(x+(\sqrt{2}-1)^2\sqrt{1-x^2})$
                                                proportional to the term depending on $x$ in $\Sigma_{11}(a,x)+|\Sigma_{12}(a,x)|$ is increasing on $[\frac{2^{5/4}}{1+\sqrt{2}},\frac{1+\sqrt{2}}{\sqrt{6}}]$ and decreasing on $[\frac{1+\sqrt{2}}{\sqrt{6}},1]$. We deduce that when $(a,x)\in[\frac{17}{3},3+2\sqrt{2}]\times[\frac{2^{5/4}}{1+\sqrt{2}},\frac{1+\sqrt{2}}{\sqrt{6}}]$ is such that $a+2(1+\sqrt{2})x\le \frac{17}{3}+2\frac{(1+\sqrt{2})^2}{\sqrt{6}}$ i.e. $\Sigma_{11}(a,x)\le\Sigma_{11}(\frac{17}3,f(\frac{17}3))$, $\Sigma(a,x)$ is dominated and % is smaller than some matrix $\left(\begin{array}{cc}1+\frac{17}3+2(1+\sqrt{2})f(\frac{17}3) & b\\b & 1+\frac{17}3+2(1+\sqrt{2})f(\frac{17}3)\end{array}\right)$ satisfying \eqref{inecov} and therefore
                                                satisfies \eqref{inecov}. Moreover, $[\frac{17}3,3+2\sqrt{2}]\ni a\mapsto \Sigma_{11}(a,f(a))+|\Sigma_{12}(a,f(a))|$ is decreasing while $[\frac{17}3,3+2\sqrt{2}]\ni a\mapsto \Sigma_{11}(a,f(a))$ is increasing. As a consequence, when $a+2(1+\sqrt{2})x> \frac{17}{3}+2\frac{(1+\sqrt{2})^2}{\sqrt{6}}$, then $\Sigma(a,x)$ is dominated if and only if $|\Sigma_{12}(a,x)|\le |\Sigma_{12}(\hat a,f(\hat a))|$ where $\hat a$ is such that $a+2(1+\sqrt{2})x=\hat a+2(1+\sqrt{2})f(\hat a)$ i.e. $\Sigma_{11}(a,x)=\Sigma_{11}(\hat a,f(\hat a))$. Indeed, $\Sigma_{11}(\tilde a,f(\tilde a))<\Sigma_{11}(a,x)$ for $\tilde a<\hat a$ and $\Sigma_{11}(\tilde a,f(\tilde a))+|\Sigma_{12}(\tilde a,f(\tilde a))|<\Sigma_{11}(\hat a,f(\hat a))+|\Sigma_{12}(\hat a,f(\hat a))|$ for $\tilde a>\hat a$. For  $(a,x)=(\frac{17}3,1)$, $\hat a=5.7152$ and $|\Sigma_{12}(\hat a,f(\hat a))|=0.3972>\frac 13=|\Sigma_{12}(\frac{17}3,1)|$ so that $\Sigma(\frac{17}3,1)$ satisfies \eqref{inecov}.

                                                The function $$y\mapsto \frac{d}{dy}\left|\Sigma_{12}\left(a-y,f(a)+\frac{y}{2(1+\sqrt{2})}\right)\right|=\frac{a-y-3}{8\sqrt{1-\frac{(a-y-3)^2}8}}-(\sqrt{2}-1)^2\frac{f(a)+\frac{y}{2(1+\sqrt{2})}}{\sqrt{1-(f(a)+\frac{y}{2(1+\sqrt{2})})^2}}$$ is decreasing on $[(a-3-2\sqrt{2})\vee 2(1+\sqrt{2})(\frac{2^{5/4}}{1+\sqrt{2}}-f(a)) ,(a-\frac{17}3)\wedge 2(1+\sqrt{2})(1-f(a))]$, equal to $\frac{(9-6\sqrt{2})(a-\frac{17}3)}{8\sqrt{1-\frac{(a-3)^2}8}}>0$ for $y=0$ and is equal to $0$ when $f(a)+\frac{y}{2(1+\sqrt{2})}=g(a-y)$ with $g(\tilde a)=\left(1+\frac 8{(1+\sqrt{2})^4}\left(\frac{8}{(\tilde a-3)^2}-1\right)\right)^{-1/2}$. We deduce that for $y\le 0$ such that $(a-y,f(a)+\frac{y}{2(1+\sqrt{2})})\in [\frac{17}{3},3+2\sqrt{2}]\times[\frac{2^{5/4}}{1+\sqrt{2}},1]$, $\Sigma(a-y,f(a)+\frac{y}{2(1+\sqrt{2})})$ is dominated. The two previous domination cases ensure that for  $(a,x)\in[\frac{17}{3},3+2\sqrt{2}]\times[\frac{2^{5/4}}{1+\sqrt{2}},1]$ such that $x\le f(a)$, $\Sigma(a,x)$ is dominated and satisfies \eqref{inecov}. On the other hand, for $a\in (\frac{17}3,3+2\sqrt{2})$ and $y>0$ small enough, $\Sigma(a-y,f(a)+\frac{y}{2(1+\sqrt{2})})$ is not dominated. % Indeed, for $y>0$ small enough, $|\Sigma_{12}(a-y,f(a)+\frac{y}{2(1+\sqrt{2})})|>|\Sigma_{12}(a,f(a))|$ while clearly $\Sigma_{11}(a-y,f(a)+\frac{y}{2(1+\sqrt{2})})=\Sigma_{11}(a,f(a))$. Then $\Sigma(a-y,f(a)+\frac{y}{2(1+\sqrt{2})})$ is not smaller than $\left(\begin{array}{cc}\Sigma_{11}(\tilde a,f(\tilde a))& b\\b & \Sigma_{11}(\tilde a,f(\tilde a))\end{array}\right)$ for some $b$ such that $|b|\le|\Sigma_{12}(\tilde a,f(\tilde a))|$ when $\tilde a=a$. It is not smaller when $\tilde a<a$ since $\tilde a\mapsto \Sigma_{11}(\tilde a,f(\tilde a))$ is increasing and when $\tilde a>a$ since $a\mapsto \Sigma_{11}(a,f(a))+|\Sigma_{12}(a,f(a))|$ is decreasing. 
                                                In particular, for $a\in(\frac{17}3,3+2\sqrt{2})$, $\Sigma(a,g(a))$ is not dominated (note that $g$ is increasing on $[\frac{17}3,3+2\sqrt{2}]$ with $g(\frac{17}3)=\frac{1+\sqrt{2}}{\sqrt{6}}$ and $g(3+2\sqrt{2})=1$). % The couples $(a,x)$ such that $\Sigma(a,x)$ is not dominated are represented in Figure \ref{figax} \textcolor{red}{exgacc.sce} with the upper boundary obtained numerically by searching for $y>0$ such that $|\Sigma_{12}(a-y,f(a)+\frac{y}{2(1+\sqrt{2})})|=|\Sigma_{12}(a,f(a))|$.
% The maximal absolute extradiagonal coefficient of $\Sigma(a,x)$ for fixed diagonal coefficient is attained along the curve $(a,g(a)$ where $g(a)=\left(1+\frac 8{(1+\sqrt{2})^4}\left(\frac{8}{(a-3)^2}-1\right)\right)^{-1/2}$.

\end{Example}

                                                                                                                                                                                                                                         \subsection{Comparison of the sufficient conditions with \eqref{condconvjp}}
We now compare the sufficient conditions \eqref{inecov} and \eqref{correl} to the generalization of \eqref{condconvjp} to $n\ge 2$. 
\begin{Proposition}
  \label{proplieno}
      Assume that $\Sigma=\sigma\sigma^*,\Sigma_1=\sigma_1\sigma_1^*,\cdots,\Sigma_n=\sigma_n\sigma_n^*$ for $\sigma,\sigma_1,\cdots,\sigma_n\in\R^{d\times q}$. Then \begin{equation}
 \exists O_1,\cdots,O_n\in{\cal O}(q)\mbox{ such that }\sigma\sigma^*\le \left(\sum_{i=1}^np_i\sigma_iO_i\right)\left(\sum_{i=1}^np_i\sigma_iO_i\right)^*\label{condsqrt}
\end{equation}
implies \eqref{inecov} and, when $q\ge nd$, the converse implication holds. Moreover, when $q\ge d$, \eqref{correl} implies \eqref{condsqrt}. % implies \eqref{inecov} while, when $q\ge d$, \eqref{correl} implies \eqref{condsqrt}. 
In particular, \eqref{correl} implies 
\begin{equation}\exists O_1,\cdots,O_n\in{\cal O}(d)\mbox{ such that }\sigma\sigma^*\le \left(\sum_{i=1}^np_i\Sigma_i^{1/2}O_i\right)\left(\sum_{i=1}^np_i\Sigma_i^{1/2}O_i\right)^*.\label{condsqrtsqrt}\end{equation}
\end{Proposition}
\begin{Remark}One has
  $$\left(\sum_{i=1}^np_i\sigma_iO_i\right)\left(\sum_{i=1}^np_i\sigma_iO_i\right)^*=\left(p_1\sigma_1+\sum_{i=2}^np_i\sigma_iO_iO_1^*\right)\left(p_1\sigma_1+\sum_{i=2}^np_i\sigma_iO_iO_1^*\right)^*,$$
  with $O_iO_1^*\in{\cal O}(q)$ when $O_1,O_i\in{\cal O}(q)$. As a consequence, we can choose $O_1=I_q$ in \eqref{condsqrt} without making the condition stronger. In particular, for $n=2$, \eqref{condsqrt} is equivalent to
  $$\exists O\in{\cal O}(q),\;\sigma\sigma^*\le \left(p_1\sigma_1+(1-p_1)\sigma_2 O\right)\left(p_1\sigma_1+(1-p_1)\sigma_2 O\right)^*,$$
  which, when $(\sigma,\sigma_1,\sigma_2)$ is replaced by $(\sigma(p_1x+(1-p_1)y),\sigma(x),\sigma(y))$, writes \eqref{condconvjp}.
\end{Remark}We deduce that under equivalence of \eqref{correl} and \eqref{compccgg} and in particular in the four settings given at the end of Theorem \ref{proprinc},  then \eqref{condsqrtsqrt} and \eqref{compccgg} also are equivalent and, when $q\ge d$, \eqref{condsqrt} and \eqref{compccgg} also are equivalent. Since for $U_i\in\R^{d\times nd}$ such that $U_iU_i^*=I_d$, $(\Sigma^{1/2}_i U_i)(\Sigma^{1/2}_i U_i)^*=\Sigma_i$ and $(U_iO_i)(U_iO_i)^*=I_d$ when $O_i\in{\cal O}(nd)$ , we also deduce the following corollary.
 \begin{Corollary} The condition \eqref{inecov} is equivalent to 
  \begin{align}
   \exists U_1,\cdots,U_n\in\R^{d\times nd}&\mbox{ such that }U_iU_i^*=I_d\mbox{ for }i\in\{1,\cdots,n\}\notag\\&\mbox{ and  }\Sigma\le \left(\sum_{i=1}^n p_i\Sigma_i^{1/2}U_i\right)\left(\sum_{i=1}^n p_i\Sigma_i^{1/2}U_i\right)^*.\label{inecovnd}
  \end{align}
\end{Corollary}

\noindent{\bf Proof of Proposition \ref{proplieno}.}
The condition \eqref{condsqrt} implies \eqref{inecov} for $\Gamma=BB^*$ where $B\in\R^{nd\times q}$ is defined by $B_{(i-1)d+1:id,1:q}=\sigma_i O_i$ for $i\in\{1,\cdots,n\}$ so that $(p_1I_d,\cdots,p_nI_d)B=\sum_{i=1}^np_i\sigma_iO_i$. When $q\ge nd$, to prove the converse implication, we suppose \eqref{inecov} and  for $i\in\{1,\cdots,n\}$, we denote by
$\Theta_i\in\R^{d\times q}$ the matrix such that $(\Theta_i)_{1:d,1:nd}=\Gamma^{1/2}_{(i-1)d+1:id,1:nd}$ and with all other entries equal to $0$.
Then, for $i,j\in\{1,\cdots,n\}$,  $$\Theta_i\Theta_j^*=\Gamma^{1/2}_{(i-1)d+1:id,1:nd}(\Gamma^{1/2}_{(j-1)d+1:jd,1:nd})^*=\Gamma_{(i-1)d+1:id,(j-1)d+1:jd}.$$ In particular, $\Theta_i\Theta_i^*=\Sigma_i=\sigma_i\sigma_i^*,$ so that, according to Lemma A.1 \cite{JourPagVolt22}, there exists $O_i\in{\cal O}(q)$ such that $\Theta_i=\sigma_iO_i$. Moreover,
$$\Sigma\le (p_1I_d,\cdots,p_nI_d)\Gamma(p_1I_d,\cdots,p_nI_d)^*=\left(\sum_{i=1}^np_i\Theta_i\right)\left(\sum_{i=1}^np_i\Theta_i\right)^*=\left(\sum_{i=1}^np_i\sigma_iO_i\right)\left(\sum_{i=1}^np_i\sigma_iO_i\right)^*.$$

% On the other hand, under \eqref{inecov}  and  for $i\in\{1,\cdots,n\}$, $\Gamma^{1/2}_{(i-1)d+1:id,1:nd}=\Sigma_i^{1/2}P_i$ and $\sigma_i=\Sigma_i^{1/2}Q_i$ where $P_i\in\R^{d\times nd}$ and $Q_i\in\R^{d\times q}$ have ${\rm dim}({\rm Im}(\Sigma_i^{1/2}))$ columns which form an orthonormal basis of ${\rm Im}(\Sigma_i^{1/2})$ and their other columns vanish.

Let us now assume \eqref{correl}. For $i\in\{1,\cdots,n\}$ and $D_i={\rm diag}\left(\sqrt{(M\Sigma_iM^*)_{11}},\cdots,\sqrt{(M\Sigma_iM^*)_{dd}}\right)$, $\Sigma_i=M^{-1}D_iCD_i(M^*)^{-1}=\left(M^{-1}D_iC^{1/2}\right)\left(M^{-1}D_iC^{1/2}\right)^*$ while for $D=\sum_{i=1}^n p_i D_i$,
  \begin{equation}
   \Sigma\le M^{-1}DCD(M^*)^{-1}=\left(\sum_{i=1}^np_iM^{-1}D_iC^{1/2}\right)\left(\sum_{i=1}^np_iM^{-1}D_iC^{1/2}\right)^*.\label{correlbis}
 \end{equation}
According to Lemma A.1 \cite{JourPagVolt22}, for $i\in\{1,\cdots,n\}$, $\sigma_i=\Sigma_i^{1/2}P_i$ and  $M^{-1}D_iC^{1/2}=\Sigma_i^{1/2}Q_i$ for some matrices $P_i\in\R^{d\times q}$ and $Q_i\in\R^{d\times d}$ such that $P_iP_i^*$ and $Q_iQ_i^*$ are the orthogonal projection on ${\rm Im}(\Sigma_i)$. When $q\ge d$, we set $\tilde Q_i\in\R^{d\times q}$ to be the matrix such that $(\tilde Q_i)_{1:d,1:d}=Q_i$ and with all other entries equal to $0$. Then $\tilde Q_i\tilde Q_j^*=Q_iQ_j^*$ for all $i,j\in\{1,\cdots,n\}$. The equality for $j=i$ implies that $\tilde Q_i\tilde Q_i^*=P_iP_i^*$ so that, by Lemma A.1 \cite{JourPagVolt22}, $\tilde Q_i=P_iO_i$ for some $O_i\in{\cal O}(q)$. Hence for $i,j\in\{1,\cdots,n\}$,
\begin{align*}
   \sigma_iO_i(\sigma_j O_j)^*&=\Sigma_i^{1/2}P_iO_i(P_jO_j)^*\Sigma_j^{1/2}=\Sigma_i^{1/2}\tilde Q_i\tilde Q_j^*\Sigma_j^{1/2}=\Sigma_i^{1/2}Q_iQ_j^*\Sigma_j^{1/2}\\&=M^{-1}D_iC^{1/2}(M^{-1}D_jC^{1/2})^*.
\end{align*}
Hence \eqref{correlbis} writes
$$\sigma\sigma^*\le \left(\sum_{i=1}^np_i\sigma_iO_i\right)\left(\sum_{i=1}^np_i\sigma_iO_i\right)^*.$$
\hfill$\Box$

Since \eqref{condconvjp} could accomodate a $\R^q$-valued random noise $Z$ with radial distribution in place of $G\sim{\cal N}_q(0,I_q)$, let us finally consider such a more general noise.

\begin{Proposition}\label{propradno}
  Let $Z$ be some integrable $\R^q$-valued random vector with radial distribution, $n\ge 2$, $\sigma_1,\cdots,\sigma_n\in\R^{d\times q}$ and $(p_1,\cdots,p_n)\in(0,1)^n$ be such that $\sum_{i=1}^n p_i=1$. Then 
  $$\mbox{\eqref{condsqrt}}\Longrightarrow{\cal L}(\sigma Z)\le_{cx}\sum_{i=1}^np_i{\cal L}(\sigma_i Z)\Longrightarrow\forall \xi\in\R^d,\;\sqrt{\xi^*\sigma\sigma^*\xi}\le \sum_{i=1}^n p_i\sqrt{\xi^*\sigma_i\sigma^*_i\xi}.$$
  Moreover,  
  $$\forall i\in\{1,\cdots,n\},\;\sigma_i\sigma_i^*\le \sigma\sigma^*\Longrightarrow\sum_{i=1}^np_i{\cal L}(\sigma_i Z)\le_{cx}{\cal L}(\sigma Z),$$
with equivalence if $Z$ has finite moments of all orders,\end{Proposition}
\begin{Remark}\label{remrad}
   For $Z=(Z_1,\cdots,Z_q)$ an integrable $\R^q$-valued random vector with radial distribution such that $\P(|Z|>0)>0$ and $\sigma,\tilde\sigma\in\R^{d\times q}$, one has ${\cal L}(\sigma Z)\le_{cx}{\cal L}(\tilde\sigma Z)\Leftrightarrow \sigma\sigma^*\le \tilde\sigma\tilde\sigma^*$. This was stated in \cite[Lemma 3.2 (iii)]{JP} under  square integrability of $Z$ for the necessary condition. Proposition \ref{propradno} applied with $\sigma_1=\cdots=\sigma_n=\tilde \sigma$ ensures that this additional assumption can be relaxed.
\end{Remark}

\noindent{\bf Proof.}
Under \eqref{condsqrt}, by \cite[Lemma 3.2 (iii)]{JP}, ${\cal L}(\sigma Z)\le_{cx} {\cal L}(\sum_{i=1}^np_i\sigma_iO_i Z)$. Since ${\cal L}(\sum_{i=1}^np_i\sigma_iO_i Z)\le_{cx} \sum_{i=1}^np_i{\cal L}(\sigma_iO_i Z)$ and, by the radiality of ${\cal L}(Z)$, ${\cal L}(\sigma_iO_i Z)={\cal L}(\sigma_iZ)$, this implies that ${\cal L}(\sigma Z)\le_{cx}\sum_{i=1}^np_i{\cal L}(\sigma_i Z)$.

Let $e_1,\cdots,e_q$ denote the canonical basis of $\R^q$. For $z\in\R^q\setminus\{0\}$ and $O\in{\cal O}(q)$ such that $O\frac{z}{|z|}$ is equal to $e_1$, by the radiality of ${\cal L}(Z)$,
$z^*Z$ has the same distribution as $|z|\left(\frac{z}{|z|}\right)^*O^*Z=|z|Z_1$. Moreover, $\E[|Z_1|]=\frac{1}{q}\sum_{i=1}^q\E[|e_i^*Z|]=\frac{1}{q}\E[\sum_{i=1}^q|Z_i|]>0$.
As a consequence, for $\xi\in\R^d$ and $\theta\in\R^{d\times q}$, $\E[|\xi^*\theta Z|]=|\theta^*\xi|\E[|Z_1|]$, and, for the convex function $\R^d\ni x\mapsto|\xi^*x|$, ${\cal L}(\sigma Z)\le_{cx}\sum_{i=1}^np_i{\cal L}(\sigma_i Z)$ implies $\sqrt{\xi^*\sigma\sigma^*\xi}\le \sum_{i=1}^n p_i\sqrt{\xi^*\sigma_i\sigma^*_i\xi}$.

When $\sigma_i\sigma_i^*\le \sigma\sigma^*$ for $i\in\{1,\cdots,n\}$, then, by Remark \ref{remrad}, ${\cal L}(\sigma_i Z)\le_{cx}{\cal L}(\sigma Z)$ for $i\in\{1,\cdots,n\}$ and $\sum_{i=1}^np_i{\cal L}(\sigma_i Z)\le_{cx}\sum_{i=1}^np_i{\cal L}(\sigma Z)={\cal L}(\sigma Z)$. When $Z$ has finite moments of all orders $k\ge 1$, for the convex function $\varphi(x)=|\xi^*x|^k$ with $\xi\in\R^d$, $\sum_{i=1}^np_i{\cal L}(\sigma_i Z)\le_{cx}{\cal L}(\sigma Z)$
implies that $(\E[|Z_1|^k])^{1/k}\left(\sum_{i=1}^n|\xi^*\sigma_i\sigma_i^*\xi|^{k/2}\right)^{1/k}\le (\E[|Z_1|^k])^{1/k}\sqrt{\xi^*\sigma\sigma^*\xi}$, from which we deduce in the limit $k\to\infty$ that $\max_{i\in\{1,\cdots,n\}}\sqrt{\xi^*\sigma_i\sigma_i^*\xi}\le \sqrt{\xi^*\sigma\sigma^*\xi}$.

\hfill$\Box$

                                                                                 \section{The case with $n=2$  Gaussian distributions in the mixture}                                                Let us now consider the particular case $n=2$ which was our initial motivation as seen in the introduction. In view of Lemma \ref{lemcovfsq}, this case also is of particular interest in order to check \eqref{inegsqrt} even when $n\ge 3$. By Remark \ref{remskew}, when $n=2$, \eqref{inegsqrt} also writes
\begin{equation}
   \forall \alpha>0,\;\Sigma\le p_1\Sigma_1+(1-p_1)\Sigma_2+ p_1(1-p_1)\left(\alpha\Sigma_1+\frac 1\alpha\Sigma_2\right)% p_1\left(p_1+\alpha (1-p_1)\right)\Sigma_1+(1-p_1)\left((1-p_1)+\frac{p_1}{\alpha}\right)\Sigma_2
   .\label{inegsqrtbis}
\end{equation} On the other hand, we can rewrite \eqref{inecov} in the following way.
\begin{Lemma}\label{leminecovn2}
  When $n=2$, \eqref{inecov} is equivalent to
  \begin{align}&\Sigma\le p_1^2\Sigma_1+(1-p_1)^2\Sigma_2+p_1(1-p_1)(\Theta+\Theta^*)\notag\\
   \mbox{ for some }&\Theta\in\R^{d\times d}\mbox{ such that }\forall x,y\in\R^d,\;(x^*\Theta y)^2\le x^*\Sigma_1x\times y^*\Sigma_2y.\label{inecovn2}
\end{align}
\end{Lemma}
Since any symmetric matrix $\Gamma\in\R^{2d\times 2d}$ such that $\Gamma_{1:d,1:d}=\Sigma_1$ and $\Gamma_{d+1:2d,d+1:2d}=\Sigma_2$ writes $
    \Gamma=\left(\begin{array}{cc} \Sigma_1 & \Theta
   \\ \Theta^* & \Sigma_2
                 \end{array}\right)$ for some $\Theta\in\R^{d\times d}$, this equivalence is an easy consequence of the following lemma.\begin{Lemma}\label{lem+n2}
  For $\Sigma_1,\Sigma_2\in{\cal S}_+(d)$ and $\Theta\in\R^{d\times d}$, one has \begin{align*}
   \left(\begin{array}{cc} \Sigma_1 & \Theta
   \\ \Theta^* & \Sigma_2
                 \end{array}\right)\in{\cal S}_+(2d)&\Longleftrightarrow\Theta=S_1S_2^*\mbox{ with }S_1,S_2\in\R^{d\times d}\mbox{ such that }S_1S_1^*\le\Sigma_1,\;S_2S_2^*\le\Sigma_2\\&\Longleftrightarrow \forall x,y\in\R^d,\;(x^*\Theta y)^2\le x^*\Sigma_1x\times y^*\Sigma_2y.
  \end{align*}
\end{Lemma}
\begin{Remark}
   One can choose either $S_1=\Sigma_1^{1/2}$ or $S_2=\Sigma_2^{1/2}$ in the second equivalent assertion.
\end{Remark}
\noindent{\bf Proof of Lemma \ref{lem+n2}.} Let $\Gamma=\left(\begin{array}{cc} \Sigma_1 & \Theta
   \\ \Theta^* & \Sigma_2
                                                              \end{array}\right)$.

                                                            The equality $\Theta=S_1S_2^*$ with $S_1S_1^*\le\Sigma_1$ and $S_2S_2^*\le\Sigma_2$ combined with Cauchy-Schwarz inequality imply that $$\forall x,y\in\R^d,\;(x^*\Theta y)^2=\left((S_1^*x)^*S_2^*y\right)^2\le x^*S_1S_1^*x\times y^*S_2S_2^*y\le x^*\Sigma_1x\times y^*\Sigma_2 y.$$
Under the latter property,
$$\forall x,y\in\R^d,\;(x^*,y^*)\Gamma\left(\begin{array}{c} x
                  \\ y
                                         \end{array}\right)=x^*\Sigma_1x+2x^*\Theta y+y^*\Sigma_2 y\ge \left(\sqrt{x^*\Sigma_1x}-\sqrt{y^*\Sigma_2y}\right)^2\ge 0,$$
                                       so that $\Gamma\in{\cal S}_+(2d)$.
                                       
When $\Sigma_2$ is positive definite, the convex function $\R^d\ni y\mapsto 2x^*\Theta y+y^*\Sigma_2 y$ being minimal and equal to $-x^*\Theta \Sigma_2^{-1}\Theta^*x$ at $y=-\Sigma_2^{-1}\Theta^*x$ for fixed $x\in\R^d$, we recover that $\Gamma\in{\cal S}_+(2d)$ is equivalent to the fact that the Schur complement $\Sigma_1-\Theta\Sigma_2^{-1}\Theta^*=\Sigma_1-\Theta\Sigma_2^{-1/2}(\Theta\Sigma_2^{-1/2})^*$ belongs to ${\cal S}_+(d)$. Then $S_1=\Theta\Sigma_2^{-1/2}$ satisfies $S_1S_1^*\le \Sigma_1$ and $\Theta=S_1S_2^*$ with $S_2=\Sigma_2^{1/2}$ such that $S_2S_2^*=\Sigma_2$. Let us finally suppose that $\Gamma\in{\cal S}_+(2d)$. For each $k\in\N^*$, $\Gamma_k=\Gamma+\frac 1k I_{2d}=\left(\begin{array}{cc} \Sigma_1 +\frac 1k I_{d}& \Theta
   \\ \Theta^* & \Sigma_2+\frac 1k I_{d}
                                                                                                                                                                                                                                                                                                                                                                                                                                                                                                                                                                                                                                                                                                                         \end{array}\right)$ is positive definite with $\Sigma_2+\frac 1 k I_d$ positive definite. Therefore $\Theta=S_1^k(\Sigma_2+\frac 1 k I_d)^{1/2}$ with $S_1^k\in\R^{d\times d}$ such that $S_1^k(S_1^k)^*\le \Sigma_1+\frac 1k I_d$. Since for $i\in\{1,\cdots,d\}$ and $k\in\N^*$, $\sum_{j=1}^d (S_1^k)_{ij}^2\le (\Sigma_1)_{ii}+\frac 1k$ and $\sum_{j=1}^d ((\Sigma_2+\frac 1 k I_d)^{1/2}_{ij})^2=(\Sigma_2)_{ii}+\frac 1k$, we can extract from $(S_1^k,(\Sigma_2+\frac 1 k I_d)^{1/2})_{k\in\N^*}$ a subsequence such that all entries of the matrices converge to the entries of some limiting matrices $(S_1,S_2)$. The matrix $S_1$ is such that $S_1S_1^*\le \Sigma_1$ while the matrix $S_2$ is symmetric and such that $S_2S_2^*=\Sigma_2$ and theferore $S_2=\Sigma_2^{1/2}$. Moreover, $\Theta=S_1S_2^*$, which proves the necessary condition.\hfill$\Box$
\begin{Remark}When $\Gamma=\left(\begin{array}{cc} \Sigma_1 & \Theta
   \\ \Theta^* & \Sigma_2\end{array}\right)$ belongs to ${\cal S}_+(2d)$ with $\Sigma_2$ positive definite and $X\sim{\cal N}_{2d}(0,\Gamma)$, then the Schur complement $\Sigma_1-\Theta\Sigma_2^{-1}\Theta^*$ is the covariance matrix of the conditional law of $X_{1:d}$ given $X_{d+1:2d}$.
   
\end{Remark}
                                                                                                                                                                                                                                                                                                                                                                                                                                                                                                                                                                                                                                                                                                                   
\begin{Remark}
   If there exists $\Sigma_3\in{\cal S}_+(d)$ such that $\Sigma_3\le \Sigma_1$, $\Sigma_3\le \Sigma_2$  and $\Sigma\le p_1^2\Sigma_1+(1-p_1)^2\Sigma_2+2p_1(1-p_1)\Sigma_3$ then \eqref{inecov} holds with $\Gamma=\left(\begin{array}{cc} \Sigma_1 & \Sigma_3
   \\ \Sigma_3 & \Sigma_2
                                                                                                                                                                                                                         \end{array}\right)$ since for $S_1=S_2=\Sigma_3^{1/2}$, $S_1S_1^*=\Sigma_3\le \Sigma_1$ and $S_2S_2^*=\Sigma_3\le \Sigma_2$. The symmetric matrix $B=\frac 1{2p_1(1-p_1)}\left(\Sigma-p_1^2\Sigma_1-(1-p_1)^2\Sigma_2\right)$ writes $O{\rm diag}(\lambda_1,\cdots,\lambda_d)O^*$ for some $O\in{\cal O}(d)$. Setting $B^+=O{\rm diag}(\lambda_1^+,\cdots,\lambda_d^+)O^*$, we deduce that a sufficient condition for \eqref{inecov} to hold is $B^+\le \Sigma_1$ and $B^+\le \Sigma_2$. Note that $B\le \Sigma_1$ does not imply $B^+\le\Sigma_1$, as seen from the example $B=\left(\begin{array}{cc}1 & -1
   \\ -1 & -1
\end{array}\right)$ with  $B^+=\left(\begin{array}{cc}\frac 1{2(\sqrt{2}-1)} & -\frac 12
   \\ -\frac 12 & \frac{\sqrt{2}-1}2
\end{array}\right)$ and  $\Sigma_1=\left(\begin{array}{cc}1 & -1
   \\ -1 & 1
                                         \end{array}\right)$.

                                       On the other hand when $\Theta=S_1S_2^*$ with $S_1S_1^*=\Sigma_1$ and $S_2S_2^*=\Sigma_2$, then $\left(\begin{array}{cc} \Sigma_1 & \Theta
   \\ \Theta^* & \Sigma_2
                 \end{array}\right)=\left(\begin{array}{c} S_1 \\ S_2
                 \end{array}\right)\left(\begin{array}{c} S_1 \\ S_2
                 \end{array}\right)^*$ is the covariance matrix of $\left(\begin{array}{c} S_1 \\ S_2
                 \end{array}\right)G$ where $G\sim{\cal N}_d(0,I_d)$. Note that when $\Sigma_2$ is  positive definite then the Schur complement $\Sigma_1-\Theta\Sigma_2^{-1}\Theta^*$ vanishes if and only if $\Theta=S_1S_2^*$ with $S_1S_1^*=\Sigma_1$ and $S_2S_2^*=\Sigma_2$.
               When $\Sigma_1$ (resp. $\Sigma_2$) is positive definite, the choice $S_1=\Sigma_1^{1/2}$ and $S_2=\Sigma_1^{-1/2}(\Sigma_1^{1/2}\Sigma_2 \Sigma_1^{1/2})^{1/2}$ (resp. $S_2=\Sigma_2^{1/2}$ and $S_1=\Sigma_2^{-1/2}(\Sigma_2^{1/2}\Sigma_1 \Sigma_2^{1/2})^{1/2}$) is such that ${\cal N}_{2d}\left(0,\left(\begin{array}{c} S_1 \\ S_2
                 \end{array}\right)\left(\begin{array}{c} S_1 \\ S_2
                 \end{array}\right)^*\right)$ is the quadratic Wasserstein optimal coupling between ${\cal N}_d(0,\Sigma_1)$ and ${\cal N}_d(0,\Sigma_2)$ (see \cite{Dowlan,Olpuk,GiSho}).
Since for $\varphi:\R^d\to\R$ continuously differentiable with a Lipschitz continuous gradient and $x,y\in\R^d$,
  \begin{align*}
   &p_1\varphi(x)+(1-p_1)\varphi(y)-\varphi(p_1 x+(1-p_1)y)\\&=p_1(1-p_1)\int_{q=0}^1(x-y).\left(\nabla\varphi((1-q(1-p_1))x+q(1-p_1)y)-\nabla\varphi((1-q p_1)y+q p_1 x))\right)dq\\&\le \frac{p_1(1-p_1)}{2}{\rm Lip}(\nabla \varphi)|x-y|^2,
  \end{align*}
  the choice $(X,Y)$ distributed according to the optimal quadratic Wasserstein coupling between ${\cal L}(X)$ and ${\cal L}(Y)$ is very sensible in order to minimize the loss in the inequality ${\cal L}(p_1X+(1-p_1)Y)\le_{cx} p_1{\cal L}(X)+(1-p_1){\cal L}(Y)$, which is the second step in the derivation of \eqref{inecov}$\Rightarrow$\eqref{compccgg} in the proof of Theorem \ref{proprinc}.

  There also exist $L^{\Sigma_1},L^{\Sigma_2}\in\R^{d\times d}$ lower triangular such that $L^{\Sigma_1}{L^{\Sigma_1}}^*={\Sigma_1}$ and $L^{\Sigma_2} {L^{\Sigma_2}}^*={\Sigma_2}$ (the existence of $L^{\Sigma_1}$ and $L^{\Sigma_2}$ is guaranteed by the Cholesky decomposition when ${\Sigma_1}$ and ${\Sigma_2}$ are positive definite).
When the diagonal entries of $L^{\Sigma_1}$ and $L^{\Sigma_2}$ are positive, ${\cal N}_{2d}\left(0,\left(\begin{array}{c} L^{\Sigma_1} \\ L^{\Sigma_2}
                 \end{array}\right)\left(\begin{array}{c} L^{\Sigma_1} \\ L^{\Sigma_2}
                 \end{array}\right)^*\right)$ is the Knothe's rearrangement coupling between ${\cal N}_d(0,{\Sigma_1})$ and ${\cal N}_d(0,{\Sigma_2})$. Indeed, then for $G\sim{\cal N}_d(0,I_d)$,  $\left(\begin{array}{c}(L^{\Sigma_1} G)_1
         \\(L^{\Sigma_2} G)_1
      \end{array}\right)\sim{\cal N}_2\left(0,\left(\begin{array}{cc} (L^{\Sigma_1}_{11})^2  &L^{\Sigma_1}_{11}L^{\Sigma_2}_{11}\\L^{\Sigma_1}_{11}L^{\Sigma_2}_{11} &(L^{\Sigma_2}_{11})^2
      \end{array}\right)\right)$ is comonotonically coupled ($(L^{\Sigma_2} G)_1=\frac{L^{\Sigma_2}_{11}}{L^{\Sigma_1}_{11}}(L^{\Sigma_1}G)_1$). Next, for $i\in\{1,\cdots,d-1\}$, the conditional law $${\cal N}_2\left(\left(\begin{array}{c}
        L^{\Sigma_1}_{i+1,1:i}G_{1:i} \\L^{\Sigma_2}_{i+1,1:i}G_{1:i}
      \end{array}\right),\left(\begin{array}{cc} (L^{\Sigma_1}_{i+1 i+1})^2  &L^{\Sigma_1}_{i+1 i+1}L^{\Sigma_2}_{i+1 i+1}\\L^{\Sigma_1}_{i+1 i+1}L^{\Sigma_2}_{i+1 i+1} &(L^{\Sigma_2}_{i+1 i+1})^2
      \end{array}\right)\right)$$ of $\left(\begin{array}{c}(L^{\Sigma_1} G)_{i+1}
         \\(L^{\Sigma_2} G)_{i+1}
      \end{array}\right)$ given $((L^{\Sigma_1}G)_{1:i},(L^{\Sigma_2} G)_{1:i})$ also is a comonotonic coupling.
 
\item \end{Remark}

% \begin{Remark}
%    Note that for $\Theta\in\R^{d\times d}$, $\Gamma=\left(\begin{array}{cc}(\Theta\Theta^*)^{1/2} & \Theta\\ \Theta^* & (\Theta^*\Theta)^{1/2}\end{array}\right)\in{\cal S}_+(2d)$ so that $\forall x,y\in\R^d,\;(x^*\Theta y)^2\le x^*(\Theta\Theta^*)^{1/2}x\times y^*(\Theta^*\Theta)^{1/2}y$. Indeed, by Lemma A.1 \cite{JourPagVolt22}, $\Theta=U(\Theta^*\Theta)^{1/2}$ for some $U\in{\cal O}(d)$ and $$\Gamma=\left(\begin{array}{cc}U & 0\\ 0& I_d\end{array}\right)\left(\begin{array}{cc}(\Theta^*\Theta)^{1/2} & (\Theta^*\Theta)^{1/2}\\ (\Theta^*\Theta)^{1/2}& (\Theta^*\Theta)^{1/2}\end{array}\right)\left(\begin{array}{cc}U & 0\\ 0 & I_d\end{array}\right)^*$$ since $U(\Theta^*\Theta)^{1/2}U^*=(\Theta\Theta^*)^{1/2}$ as the square of the symmetric matrix $U(\Theta^*\Theta)^{1/2}U^*$ is equal to $U(\Theta^*\Theta)^{1/2}(\Theta^*\Theta)^{1/2} U^*=\Theta\Theta^*$. Note that the rank of $\Gamma$ is equal to that of $(\Theta^*\Theta)^{1/2}$.
% \end{Remark}

\small
\bibliographystyle{plain}\bibliography{propsdeconvex}

\end{document}